\documentclass[11pt]{amsart}
\usepackage{amsmath,amsfonts,amssymb,mathrsfs}
\usepackage{colortbl}
\definecolor{black}{rgb}{0.0, 0.0, 0.0}
\definecolor{red}{rgb}{1.0, 0.5, 0.5}

\title[   ]{Global well-posedness of the spatially homogeneous Kolmogorov-Vicsek model as a gradient flow}

\author[Figalli]{Alessio Figalli}
\address[Alessio Figalli]{\newline Department of Mathematics, \newline The University of Texas at Austin, 
\newline Speedway 2515 Stop C1200, \newline Austin, TX 78712, USA}
\email{figalli@math.utexas.edu}

\author[Kang]{Moon-Jin Kang}
\address[Moon-Jin Kang]{\newline Department of Mathematics, \newline The University of Texas at Austin, \newline Speedway 2515 Stop C1200, \newline Austin, TX 78712, USA}
\email{moonjinkang@math.utexas.edu}

\author[Morales]{Javier Morales}
\address[Javier Morales]{\newline Department of Mathematics, \newline The University of Texas at Austin, \newline Speedway 2515 Stop C1200, \newline Austin, TX 78712, USA}
\email{jmorales@math.utexas.edu}

\newtheorem{theorem}{Theorem}[section]
\newtheorem{lemma}{Lemma}[section]

\newtheorem{proposition}{Proposition}[section]
\newtheorem{remark}{Remark}[section]

\newcommand{\bbr}{\mathbb R}
\newcommand{\bbs}{\mathbb S}

\newcommand{\bbp} {\mathbb P}

\numberwithin{figure}{section}
\newcommand{\beq}{\begin{equation}}
\newcommand{\eeq}{\end{equation}}
\newcommand{\bsp}{\begin{split}}
\newcommand{\esp}{\end{split}}

\def\eps{\varepsilon }

\newcommand\adots{\mathinner{\mkern2mu\raise1pt\hbox{.}
\mkern3mu\raise4pt\hbox{.}\mkern1mu\raise7pt\hbox{.}}}

\def\charf {\mbox{{\text 1}\kern-.30em {\text l}}}
\begin{document}
%%%%%%%%%%%%%%%%

\date{\today}


\thanks{\textbf{Acknowledgment.} A. Figalli is supported by the NSF Grants DMS-1262411 and 
DMS-1361122.
M.-J. Kang is supported by Basic Science Research Program through the National Research Foundation of Korea (NRF-2013R1A6A3A03020506), and by AMS-Simons Travel Grant. 
}

\begin{abstract}
We consider the so-called spatially homogenous Kolmogorov-Vicsek model,
a non-linear Fokker-Planck equation of self-driven stochastic particles with orientation interaction under the space-homogeneity. We prove the global existence and uniqueness of weak solutions to the equation. We also show that weak solutions exponentially converge to a steady state, which has the form of the Fisher-von Mises distribution.
\end{abstract}
\maketitle \centerline{\date}

%\tableofcontents
\section{Introduction}
In this paper we study the dynamics of the probability density function $\rho(t,\omega)$, as one-particle distribution at time $t$ with direction $\omega\in\bbs^{d-1}$ (unit sphere of $\bbr^{d}$), which satisfies the system
\begin{align}
\begin{aligned}\label{main} 
&\partial_{t}\rho= \Delta_{\omega}\rho-\nabla_{\omega}\cdot\Big(\rho\,\bbp_{\omega^{\perp}}\Omega_{\rho}\Big),\\
&\Omega_{\rho}=\frac{J_{\rho}}{|J_{\rho}|}, \quad J_{\rho}=\int_{\bbs^{d-1}}\omega\,\rho\, d\omega.
\end{aligned}
\end{align}
Here the operators $\nabla_{\omega}$ and $\Delta_{\omega}$ denote the gradient and the Laplace-Beltrami operator on the sphere $\bbs^{d-1}$, respectively. The term $\bbp_{\omega^{\perp}}\Omega$ denotes the projection of the vector $\Omega$ onto the normal plane to $\omega$, describing the mean-field force that governs the orientational interaction of self-driven particles by aligning them with the direction $\Omega$ determined by the flux $J$.
Notice that $\Omega$ is not defined when $J=0$, and this singularity in the vector field is one of the main difficulties
when studying the system \eqref{main}.\\

The equation \eqref{main} is the spatially homogeneous version of the kinetic Kolmogorov-Vicsek model, which was formally derived by Degond and Motsch \cite{D-M} as a mean-field limit of the discrete Vicsek model \cite{A-H, C-K-J-R-F, G-C, Vicsek} with stochastic dynamics. Recently, the stochastic Vicsek model has received extensive attention in the mathematical topics such as the mean-field limit, hydrodynamic limit, and phase transition. Bolley, Ca$\tilde{\mbox{n}}$izo and Carrillo \cite{B-C-C} have rigorously justified the mean-field limit 
when the unit vector $\Omega$ in the force term of \eqref{main}  is replaced by a more regular vector-field,
and Degond, Frouvelle and Liu \cite{D-F-L-1} provided a complete and rigorous description of phase transitions when $\Omega$
is replaced by $\nu(|J|)\Omega$, and there is a noise intensity $\tau(|J|)$ in front of $\Delta_{\omega}\rho$, where the functions $\nu$ and $\tau$ satisfy 
\[
|J|\mapsto \frac{\nu(|J|)}{|J|}\quad\mbox{and}\quad |J|\mapsto \tau(|J|)\quad\mbox{are Lipschitz and bounded.}
\]
Indeed, this modification leads to the appearance of phase transitions such as the number and nature of equilibria, stability, convergence rate, phase diagram and hysteresis, which depend on the ratio between $\nu$ and $\tau$. It is important to observe that the assumptions of $\nu$ remove the singularity of $\Omega$ because $\nu(|J|)\Omega\to 0$ as $|J|\to 0$. 
This phase transition problem has been studied as well in \cite{A-H, C-K-J-R-F, D-F-L-1, D-F-L-2, F-L, G-C}.
%Whereas in our case \eqref{main}, no phase transition and $\Omega$ is not defined at $|J|= 0$.
Concerning studies on hydrodynamic descriptions of kinetic Vicsek model we refer to \cite{D-F-L-1, D-F-L-2, D-M, D-M-2, D-Y, F},
see also \cite{Bo-Ca, D-D-M, H-J-K} for other related studies.

For the well-posedness of the kinetic Kolmogorov-Vicsek model, Frouvelle and Liu \cite{F-L} have shown the well-posedness in the spatially homogeneous case with the ``regular'' force field $\bbp_{\omega^{\perp}}J$ instead of $\bbp_{\omega^{\perp}}\Omega$. Moreover they have provided the convergence rates towards equilibria by using the Onsager free energy functional and Lasalle's invariance principle, and their results have been applied in \cite{D-F-L-1}.   
On the other hand, Gamba and Kang \cite{Ga-Ka} recently proved the existence of weak solutions to the kinetic Kolmogorov-Vicsek model with the singular force field $\bbp_{\omega^{\perp}}\Omega$ under the a priori assumption of $|J|>0$, without handling the stability issues, whose difficulty is mainly coming from the facts that the momentum is not conserved and no dissipative energy functional. As a study for its numerical scheme, we refer to \cite{G-H-M}. \\

The purpose of this paper is to present the global well-posedness and large time behavior of weak solutions to the spatially homogeneous problem \eqref{main}. In order to prevent the singularity of $\Omega_{\rho}$, we shall consider initial probability densities $\rho_0$ satisfying $|J_{\rho_0}|>0$. Nonetheless, since the momentum $J$ is not conserved, the condition $|J_{\rho_0}|>0$ may not immediately ensure that $|J_{\rho}|>0$ for all time.
As we shall see, a formal computation actually does show that $|J_{\rho(t)}|\geq |J_{\rho_0}|e^{-2(d-1)t}$
(see Lemma \ref{lem-moment}).
However, since it does not seem obvious how to justify this estimate, we shall rather argue by approximation.
More precisely,
we first regularize the equation \eqref{main} by adding a small constant $\eps>0$ to the denominator of
$\Omega_{\rho}$. This allows us to look at \eqref{main} as the gradient flow with respect to Wasserstein distance of a $\eps$-perturbed free energy functional,
and we will be able to prove the well-posedness of the regularized equation using the time-discrete scheme by Jordan, Kindeleherer, and Otto \cite{JKO}. Finally, using a compactness argument, we will obtain the global well-posedness of \eqref{main}.

For the large time behavior, we observe that,
as a consequence of \eqref{formula-2}, the system \eqref{main} can be written as the nonlinear Fokker-Planck equation:
\beq\label{main-0}
\partial_{t}\rho= \Delta_{\omega}\rho-\nabla_{\omega}\cdot\Big(\rho\,\nabla_{\omega}(\omega\cdot\Omega_{\rho})\Big).
\eeq
We can easily see that the equilibrium states of \eqref{main-0} have the form of the Fisher-von Mises distribution: for any given $\Omega\in\bbs^{d-1}$,
these are given by
\[
M_{\Omega}(\omega):=C_M e^{\omega\cdot\Omega},
\]
where $C_M$ is the positive constant given by
\begin{equation}
\label{eq:CM}
C_M=\frac{1}{\int_{\bbs^{d-1}}e^{\omega\cdot\Omega}\,d\omega},
\end{equation}
so that $M_{\Omega}$ is a probability density function. Notice that the normalization constant $C_M$ does not depend on $\Omega$,
and can be easily computed when $d=3$ (see Appendix). In this paper we prove that any weak solution of \eqref{main} converges
 exponentially to a stationary Fisher-von Mises distribution.\\

The paper is organized as follows.
In the next section, we briefly present some useful results and estimates in the optimal transportation theory, and then state our main results. Section 3 is devoted to the proof of existence of weak solutions. In Section 4, we prove the convergence of weak solutions towards the equilibrium in $L^1$ distance. In Section 5, we show that weak solutions are locally stable with respect to the Wasserstein distance, and as a consequence we obtain the uniqueness of the weak solution.

\section{Preliminaries and Main results}
\setcounter{equation}{0}
\subsection{Probability measures on the sphere}
Here we summarize useful results from optimal transportation theory that will be used throughout the paper.
We consider the embedded Riemannian manifold $\bbs^{d-1}\subset\mathbb{R}^{d}$ endowed with the ambient metric and geodesic distance given by 
\[
 d(x,y):=\inf\bigg\{\sqrt{\int_{0}^{1}\mid\dot{\gamma}\mid^{2}dt}~\Big|~ \gamma\in C^{1}((0,1),\bbs^{d-1}),\gamma(0)=x,\gamma(1)=y\bigg\}.
\]
We define the 2-Wasserstein distance (or transportation distance) with quadratic cost between two probability measures $\mu$
and $\nu$ as
\begin{equation}
W_{2}(\mu,\nu):=\sqrt{\inf_{\lambda\in\Lambda(\mu,\nu)}\int_{\bbs^{d-1}\times \bbs^{d-1}}d(x,y)^{2}\,d\lambda(x,y)},
\label{distance}
\end{equation}
where $\Lambda(\mu,\nu )$ denotes the set of all probability measures $\lambda$ on $\bbs^{d-1}\times \bbs^{d-1}$ with marginals $\mu$ and $\nu$,
i.e, 
\[
\pi_{1\#}\lambda = \mu,\quad \pi_{2\#}\gamma = \nu,
\]
where $\pi_1: (x,y)\mapsto x$ and $\pi_2 : (x,y)\mapsto y $ are the natural projections from $\bbs^{d-1}\times \bbs^{d-1}$ to $\bbs^{d-1}$, and $\pi_{1\#}\lambda$ denotes the push forward of $\lambda$ through $\pi_1$.\\
Whenever $\mu$ is absolutely continuous with respect to the volume measure of $\bbs^{d-1}$, it follows by McCann's Theorem \cite{MC} that there exists a unique optimal plan $\lambda_0\in\Lambda(\mu,\nu)$ which minimizes \eqref{distance},
and such a plan is induced by an optimal transport map $T: \bbs^{d-1}\rightarrow \bbs^{d-1}$, i.e., $\lambda_0
=(Id,T)_{\#}\mu$ (thus, $T_{\#}\mu=\nu$). In addition, $T$ can be written as  
\[
T(\omega)=\exp_{\omega}(\nabla \varphi(\omega)),
\]
for some $d^2/2$-convex function $\varphi :  \bbs^{d-1} \to \bbr$
(see for instance \cite[Theorem 2.33]{user}).

We shall denote by $(\mathcal{P}(\bbs^{d-1}),W_{2})$ the metric space of probability measures on the sphere endowed with the Wasserstein distance. 
We recall that $(\mathcal{P}(\bbs^{d-1}),W_{2})$ is a complete separable compact
metric space, and a sequence $\mu_{n}$ converges to $\mu$ in $W_{2}$ if and only if it converges weakly in duality with functions in $C(\bbs^{d-1})$ (see for example \cite[Theorem 3.7 and Remark 3.8]{user}). 

The following proposition provides a useful estimate on the directional derivative of the map
$\mu\mapsto W_{2}^{2}(\mu,\nu)$, which is used in Section 3. 

\begin{proposition}\label{prop-W}
 Let $\mu,\nu\in \mathcal{P}(\bbs^{d-1})$, assume that $\mu$ is absolutely continuous,
let  $X: \bbs^{d-1}\rightarrow T\bbs^{d-1}$ be a $C^{\infty}$ vector field, and define $\mu_{t}:=\exp(tX)_{\#}\mu$.
 Then we have
\[
  \limsup_{t\rightarrow0}\frac{W_{2}^{2}(\mu_{t},\nu)-W_{2}^{2}(\mu,\nu)}{t}\leq-2\int_{\bbs^{d-1}}\nabla_{\omega}\varphi (\omega)\cdot X(\omega)\,d\mu,
\]
 where $\varphi:\bbs^{d-1}\to \bbr$ is a $d^2/2$-convex function such that
 $\exp_{\omega}(\nabla_{\omega}\varphi)$ is the optimal map sending $\mu$ onto $\nu$.
 \end{proposition}
\begin{proof}
Let $\lambda_0\in\Lambda(\mu,\nu)$ be the optimal plan, i.e., $(Id,\exp_{\omega}(\nabla_{\omega} \varphi))_{\#}\mu=\lambda_0$. Since the measure
 \[
  \lambda_{t}:=\big((\exp(tX)\circ \pi_{1},\pi_{2} \big)_{\#}\lambda_0
 \]
belongs to $\Lambda(\mu_t,\nu)$, it follows by the definition of $W_2$ (see \eqref{distance}) that
\begin{align*}
\begin{aligned} 
 W_{2}^{2}(\mu_{t},\nu) &\leq \int_{\bbs^{d-1}\times \bbs^{d-1}}d(\omega,\bar\omega)^{2}\,d\lambda_t\\
 &=\int_{\bbs^{d-1}\times \bbs^{d-1}} d(\exp_{\omega}(tX),\bar\omega)^2 \,d\lambda_0\\
 &=\int_{\bbs^{d-1}} d\bigl(\exp_{\omega}(tX),\exp_{\omega}(\nabla_{\omega} \varphi)\bigr)^2 \,d\mu.
\end{aligned}
\end{align*}
We now recall that the following formula about the squared distance function (see for instance \cite[Section 1.9]{F-V}):
\begin{equation}\label{taylor}
d\bigl(\exp_{\omega}(tX(\omega)),\exp_{\omega}(\nabla_{\omega} \varphi (\omega))\bigr)^2 \leq d\bigl(\omega,\exp_{\omega}(\nabla_{\omega} \varphi (\omega))\bigr)^2 -2tX({\omega})\cdot \nabla_{\omega} \varphi (\omega)+C\, t^{2}.
  \end{equation}
Thus
\begin{align*}
\begin{aligned} 
W_{2}^{2}(\mu_{t},\nu)&\leq\int_{\bbs^{d-1}}d\bigl(\omega,\exp_{\omega}(\nabla_{\omega} \varphi (\omega))\bigr)^2 \,d\mu
-2t\int_{\bbs^{d-1}}X({\omega})\cdot \nabla_{\omega} \varphi (\omega)\,d\mu+C\,t^2\\
&=W_{2}^{2}(\mu,\nu)-2t\int_{\bbs^{d-1}}X({\omega})\cdot \nabla_{\omega} \varphi (\omega)\,d\mu+C\,t^2,
\end{aligned}
\end{align*}
and the result follows.
\end{proof}

Throughout the paper, we mainly deal with absolutely continuous measures. Hence, by abuse of notation, we will use sometimes $\rho$ to denote the absolutely continuous measure $\rho \,d\omega$ on the sphere $\bbs^{d-1}$.

\subsection{Formulas for the calculus on the sphere} 
\label{secf:formulas}
We present here some useful formulas on sphere $\bbs^{d-1}$, which are used throughout the paper.\\
Let $F:\bbs^{d-1}\to \bbr^d$ be a vector-valued function and $f:\bbs^{d-1}\to \bbr$ be scalar-valued function. Then we have the following formulas related to the integration by parts:
\beq\label{formula-0}
\int_{\bbs^{d-1}} f\,\nabla_{\omega}\cdot F \,d\omega =  -\int_{\bbs^{d-1}} F\cdot(\nabla_{\omega}f -2\omega f)\, d\omega,
\eeq
and
\begin{align}
\begin{aligned} \label{form-0}
&\int_{\bbs^{d-1}} \omega \,\nabla_{\omega} \cdot F \,d\omega = -\int_{\bbs^{d-1}} F\, d\omega,\\
&\int_{\bbs^{d-1}} \nabla_{\omega} f \,d\omega = (d-1)\int_{\bbs^{d-1}} \omega\, f \,d\omega.
\end{aligned}
\end{align}
Since $\nabla_\omega f$ is a tangent vector-field, 
it follows immediately by the definition of the projection $\bbp_{\omega^{\perp}}$ that
\begin{align}
\begin{aligned}\label{formula-1}
&\bbp_{\omega^{\perp}}\omega = 0,\\
&\bbp_{\omega^{\perp}}\nabla_{\omega} f =\nabla_{\omega} f.
\end{aligned}
\end{align}
Moreover, for any constant vector $v\in\bbr^d$ we have
\begin{align}
\begin{aligned}\label{formula-2}
&\nabla_{\omega}(\omega\cdot v) = \bbp_{\omega^{\perp}} v,\\
&\nabla_{\omega}\cdot(\bbp_{\omega^{\perp}}v) = -(d-1) \,\omega\cdot v.
\end{aligned}
\end{align}
We refer to \cite{O-T} for the derivations of the above formulas.

\subsection{Main results}
We now state our main existence, uniqueness, and convergence results.
In the sequel we shall restrict to the case $d \geq 3$ since 
we will need to use the logarithmic Sobolev inequality on the sphere (see the proof of Lemma \ref{lem-equili}).
We point out that Lemma \ref{lem-equili} is not used in the existence part,
hence our approach allows one to get existence of solutions even in the case $d=2$.

\begin{theorem} (Existence and Uniqueness) \label{thm-exist}
Assume $d \geq 3$. Let $\rho_{0}\in\mathcal{P}(\bbs^{d-1})$ be an initial probability measure satisfying 
\beq\label{ini-assume}
|J_{\rho_{0}}|>0,\quad \int_{\bbs^{d-1}}\rho_0\log\rho_0 \,d\omega<\infty.
\eeq
Then the equation \eqref{main} has a unique weak solution $\rho \in L_{loc}^{2}([0,\infty),W^{1,1}(\bbs^{d-1}))$ starting from $\rho_0$, which is weakly continuous in time, and satisfies \eqref{main} in the weak sense: for all $\varphi\in C^{\infty}(\bbs^{d-1})$ and $0\leq t<s$,
\[
 \int_{\bbs^{d-1}}\varphi\,(\rho(s)-\rho(t))\,d\omega=\int_{t}^{s}\bigg(\int_{\bbs^{d-1}}\big[\Delta_{\omega}\varphi+\nabla_{\omega}\varphi\cdot\nabla_{\omega}(\omega\cdot\Omega_{\rho(r)})\big]\,\rho(r)\,d\omega\bigg)\,dr.
\]
Moreover for all $t>0$,
\[
|J_{\rho}|^{2}\geq |J_{\rho_{0}}|^{2}e^{-2(d-1)t}.
\]
\end{theorem}

%\begin{remark}
%The assumption on dimension $d \geq 3$ is only necessary for usage of the logarithmic Sobolev inequality in the proof of Lemma \ref{lem-equili}. The Lemma \ref{lem-equili} is used to prove Theorem \ref{thm-converge}, \ref{thm-stable}, thus the uniqueness of the above theorem, but the assumption $d \geq 3$ is not necessary for the proof of the existence in Section \ref{sec-e}.
%\end{remark}

\begin{theorem} (Convergence to steady state) \label{thm-converge}
Assume $d \geq 3$. Let $\rho_{0}\in\mathcal{P}(\bbs^{d-1})$ be an initial probability measure satisfying \eqref{ini-assume}.
Then there exist a constant vector $\Omega_{\infty} \in \bbs^{d-1}$ and a constant $C>0$,
depending only on $\rho_0$ and the dimension $d$, such that
\[
 \|\rho(t) -  M_{\Omega_{\infty}}\|_{L^1(\bbs^{d-1})} \le C\Big(\int_{\bbs^{d-1}}\rho_0\log\rho_0\, d\omega+1\Big) e^{-\frac{2(d-2)}{e^{2}}t}.
\]
\end{theorem}

\begin{remark}
Notice that, since the momentum $J_{\rho(t)}$ is not conserved in time, it is not clear how to determine the vector $\Omega_{\infty}$ from the initial data $\rho_0$. 
\end{remark}

The following theorem provides a short time stability in Wasserstein distance when two initial probability measures are close to each other.
In particular it implies uniqueness of solutions.
\begin{theorem} (Stability in Wasserstein distance) \label{thm-stable}
Assume $d \geq 3$. Let $\rho_{0}, \bar\rho_0\in\mathcal{P}(\bbs^{d-1})$ be probability measures satisfying \eqref{ini-assume} and 
\[
W_{2}(\rho_{0},\bar{\rho}_{0})\leq \frac{|J_{\rho_{0}}|}{16},
\]
and let $\rho(t)$ and $\bar \rho(t)$ denote the solutions of \eqref{main} starting from $\rho_0$ and $\bar\rho_0$, respectively.
Then there exist constants $C>0$ and $\delta>0$, depending on $\rho_{0}, \bar\rho_0$, such that
\[
W_{2}(\rho(t),\bar{\rho}(t))\le e^{\lambda t}\,W_{2}(\rho_0,\bar{\rho}_0)\qquad \forall\,t<\delta,
\]
where $\lambda:=(1+2/|J_{\rho_0}|)-(d-2)$.
\end{theorem}

\section{Existence}\label{sec-e}
\setcounter{equation}{0}
In this section, we prove the existence part in Theorem \ref{thm-exist}. For this, we first regularize the equation \eqref{main} using a parameter $\eps\in (0,1)$ to prevent the singularity of $\Omega_{\rho}$, and then take $\eps\to0$ using standard compactness argument.
It is worth noticing that the existence of solutions to the regularized system could be proved also by more standard
PDE arguments. However, we prefer to use this alternative approach since it will also provide us with some useful estimates for the limiting system.

\subsection{Regularized equation}
We first regularize \eqref{main} by adding $\eps>0$ to the denominator of $\Omega_{\rho}$ as follows:
\begin{align}
\begin{aligned}\label{eqrob} 
&\partial_{t}\rho^{\eps}=\nabla_{\omega}\cdot\Big(\rho^{\eps}\,\nabla_{\omega}\,(\log\rho^{\varepsilon}-\omega\cdot\Omega^{\eps}_{\rho^{\varepsilon}})\Big),\\
&\rho^{\eps}(0)=\rho_{0},\\
&\Omega^{\varepsilon}_{\rho^{\eps}}=\frac{J_{\rho^{\eps}}}{\sqrt{|J_{\rho^{\eps}}|^{2}+\varepsilon}},\quad J_{\rho^{\eps}}=\int_{S^{d-1}}\omega\,\rho^{\eps} \,d\omega.
\end{aligned}
\end{align}
In the next subsections, we show the existence of weak solutions to the regularized equation \eqref{eqrob} as a gradient flow with respect to Wasserstein distance of the $\eps$-perturbed free energy functional 
\[
\mathcal{E}^\varepsilon(\mu):=\begin{cases}
&\int_{\bbs^{d-1}} \rho\log \rho \,d\omega -\sqrt{| J_{\rho}|^{2}+\varepsilon}\quad\mbox{if}~\mu=\rho ~d\omega \\
&+\infty,\quad \mbox{otherwise}.\\
\end{cases}
\]
Notice that since $\rho\mapsto J_{\rho}$ is continuous with respect to $W_{2}$, the functional $\mathcal{E}^{\varepsilon}$ is lower semicontinious with respect to $W_{2}$.
The next lemma provides some useful properties on derivatives of the functional $\mathcal{E}^{\eps}$.

\begin{lemma}\label{lem-E}
For a given $\rho\in\mathcal{P}(\bbs^{d-1})$, the following results hold.\\ 
(1) For any $d^2/2$-convex function $\varphi :\bbs^{d-1}\to\bbr$, the second derivative of $\mathcal{E}^{\varepsilon}$ along
the geodesic $\rho_{t} \,d\omega:=\exp(t\nabla\varphi)_{\#}\rho \,d\omega$ at $t=0$ is given by
\begin{equation}
\begin{aligned}
\label{Hessian}
&\quad\frac{d^2}{dt^{2}}\bigg|_{t=0}\mathcal{E}^{\varepsilon}(\rho_{t})=\int {\rm tr}([D^{2}\varphi]^{T}D^{2}\varphi)\,\rho \,d\omega+(d-2)\int_{\bbs^{d-1}}|\nabla\varphi|^{2}\rho \,d\omega\\
&\qquad+\int\nabla\varphi \,D^{2}(\Omega^{\varepsilon}_\rho\cdot\omega)\,\nabla\varphi\,\rho \,d\omega
- \frac{1}{\sqrt{|J_{\rho}|^{2}+\varepsilon}}\bigg(\Big|\int\nabla\varphi\,\rho \,d\omega\Big|^{2}
-\Big(\int\Omega^{\varepsilon}_{\rho}\cdot\nabla\varphi\,\rho \,d\omega\Big)^{2} \bigg).\\
\end{aligned}
\end{equation}
(2) For any smooth vector field $X: \bbs^{d-1}\rightarrow T\bbs^{d-1}$, the directional derivative of $\mathcal{E}^{\varepsilon}$ along $\mu_{t}:=\exp(tX)_{\#}\rho \,d\omega$ at $t=0$ is given by
\beq\label{direct-E}
\lim_{t\rightarrow0}\frac{\mathcal{E}^{\varepsilon}(\mu_{t})-\mathcal{E}^{\varepsilon}(\rho)}{t}=
\int_{S^{d-1}}\nabla_{\omega}(\log\rho-\omega\cdot\Omega^{\varepsilon}_\rho)\cdot X(\omega)\,\rho\,d\omega.
\eeq
(3) The slope of $\mathcal{E}^{\varepsilon}$ is given by
\beq\label{E-slope}
|\nabla\mathcal{E}^{\varepsilon}(\rho)|:=\limsup_{\bar\rho\rightarrow\rho}\frac{(\mathcal{E^{\varepsilon}}(\bar\rho)-\mathcal{E}^{\varepsilon}(\rho))_{+}}{W_{2}(\bar\rho,\rho)}
=\sqrt{\int_{\bbs^{d-1}}|\nabla_{\omega}(\log\rho-\omega\cdot\Omega^{\varepsilon}_\rho)|^{2}\rho ~d\omega}.
\eeq
\end{lemma}
\begin{proof}
Once we prove \eqref{Hessian}, since the Hessian of the map $\omega \mapsto \Omega^{\varepsilon}_\rho\cdot\omega$ has norm bounded
by $1$ and $\mbox{tr}([D^{2}\varphi]^{T}D^{2}\varphi)\ge 0$,
we get
\beq\label{see}
\frac{d^2}{dt^{2}}\bigg|_{t=0}\mathcal{E}^{\varepsilon}(\rho_{t}) \ge - \lambda \int_{\bbs^{d-1}}|\nabla\varphi|^{2}\rho \,d\omega,
\eeq
with $\lambda=(1+\eps^{-1/2})-(d-2)$.
This means that the functional $\mathcal{E}^{\varepsilon}$ is $(-\lambda)$-convex, 
and it follows by standard theory (see for instance \cite[Chapter 10]{Ambrosio-Gigli-Savare}) that \eqref{direct-E} and \eqref{E-slope}
hold.
Thus the remaining part is devoted to the proof of \eqref{Hessian}.

We begin by noticing that, since $\rho_{t} \,d\omega:=\exp(t\nabla\varphi)_{\#}\rho \,d\omega$ is a geodesic in $W_2$, 
the couple ($\rho_t,\varphi_t$) solves the following system of continuity/Hamilton-Jacobi equation in the distributional/viscosity sense
(see for instance \cite[Chapter 13]{Villani}):
\begin{align}
\begin{aligned} \label{system}
&\partial_{t}\rho_{t}+\nabla\cdot(\rho_{t}\nabla\varphi_{t})=0,\\
&\partial_{t}\varphi_t+\frac{|\nabla\varphi_{t}|^{2}}{2}=0,
\end{aligned}
\end{align}
where $\rho_0=\rho$ and $\varphi_0=\varphi$.\\
Then, using first the continuity equation above, we have
\begin{align*}
\begin{aligned}
\frac{d}{dt}\Big(\int_{\bbs^{d-1}}\rho_t\log\rho_t \,d\omega-\sqrt{|J_{\rho_t}|^{2}+\varepsilon}\Big)&=\int\log\rho_t\,\partial_{t}\rho_t \,d\omega-\frac{J_{\rho_t}}{\sqrt{|J_{\rho_t}|^{2}+\varepsilon}}\cdot\int\omega\,\partial_{t}\rho_t \,d\omega\\
&=-\int\log\rho_t\,\nabla\cdot (\rho_t\nabla\varphi_t) \,d\omega+\int\Omega^{\varepsilon}_{\rho_t}\cdot \omega\,\nabla\cdot (\rho_t\nabla\varphi_t) \,d\omega\\
&=\int\nabla\varphi_t\cdot\nabla\log\rho_t \,\rho_t \,d\omega-\int\nabla(\omega\cdot\Omega^{\varepsilon}_{\rho_t})\cdot\nabla\varphi_t\,\rho_t \,d\omega\\
&=-\int\Delta\varphi_t\,\rho_t \,d\omega -\int\nabla(\omega\cdot\Omega^{\varepsilon}_{\rho_t})\cdot\nabla\varphi_t\,\rho_t \,d\omega.
\end{aligned}
\end{align*}
Thus,
\begin{align*}
\begin{aligned}
\frac{d^2}{dt^2}\mathcal{E}^{\varepsilon}(\rho_{t})
&=-\frac{d}{dt}\int\Delta\varphi_t\,\rho_t \,d\omega 
-\frac{d}{dt}\int\nabla(\omega\cdot\Omega^{\varepsilon}_{\rho_t})\cdot\nabla\varphi_t\,\rho_t \,d\omega\\
&=: I_1 +I_2.
\end{aligned}
\end{align*}
Using \eqref{system}, we have
\begin{align*}
\begin{aligned}
I_1 &=\int \Delta\frac{|\nabla\varphi_t|^{2}}{2}\,\rho_t\,d\omega +\int\Delta\varphi_t\nabla\cdot(\nabla\varphi_t\,\rho_t)\,d\omega\\
&=\int \Delta\frac{|\nabla\varphi_t|^{2}}{2}\,\rho_t\,d\omega -\int\nabla\Delta\varphi_t\cdot\nabla\varphi_t\,\rho_t\,d\omega\\
&= \int \mbox{tr}([\nabla^{2}\varphi_t]^{T}\nabla^{2}\varphi_t)\,\rho_t \,d\omega+\int \mbox{Ric}(\nabla\varphi_t,\nabla\varphi_t)\,\rho_t \,d\omega.
\end{aligned}
\end{align*}
where in the last equality we used the Bochner formula
\[
 \Delta\frac{\mid\nabla\varphi\mid^{2}}{2}-\nabla\varphi\cdot \nabla\Delta\varphi
 =\mbox{tr}([\nabla^{2}\varphi]^{T}\nabla^{2}\varphi)+\mbox{Ric}(\nabla\varphi,\nabla\varphi).
\]
Since the Ricci curvature tensor of $\bbs^{d-1}$ is $(d-2)I_{d-1}$, we have
\[
I_1=\int \mbox{tr}([\nabla^{2}\varphi_t]^{T}\nabla^{2}\varphi_t)\,\rho_t \,d\omega+(d-2)\int |\nabla\varphi_t|^2\, \rho_t \,d\omega.
\]
For $I_2$, we use \eqref{formula-1} and \eqref{formula-2} to get
\beq\label{equ-1}
\nabla(\omega\cdot\Omega^{\varepsilon}_{\rho_t})\cdot\nabla\varphi_t = \bbp_{\omega^{\perp}}\Omega^{\varepsilon}_{\rho_t}\cdot\nabla\varphi_t =  \bbp_{\omega^{\perp}}\nabla\varphi_t\cdot\Omega^{\varepsilon}_{\rho_t} = \nabla\varphi_t\cdot\Omega^{\varepsilon}_{\rho_t},
\eeq
which yields 
\begin{align*}
\begin{aligned}
I_2 &=-\frac{d}{dt}\int\nabla\varphi_t\cdot\Omega^{\varepsilon}_{\rho_t}\,\rho_t \,d\omega\\
&=-\int\partial_t\rho_t \,\nabla\varphi_t\cdot\Omega^{\varepsilon}_{\rho_t} \,d\omega-\int\rho_t \,\nabla\partial_t\varphi_t\cdot\Omega^{\varepsilon}_{\rho_t} \,d\omega-\int\rho_t  \,\nabla\varphi_t\cdot \partial_t\Omega^{\varepsilon}_{\rho_t} \,d\omega\\
&=: I_{21} + I_{22} + I_{23}.
\end{aligned}
\end{align*}
Using \eqref{system} and \eqref{equ-1}, we have
\begin{align*}
\begin{aligned}
I_{21}&= \int  \nabla\cdot (\rho_t\nabla\varphi_t) \,\nabla\varphi_t\cdot\Omega^{\varepsilon}_{\rho_t} \,d\omega\\
&= \int \nabla\cdot (\rho_t\nabla\varphi_t) \,\nabla(\omega\cdot\Omega^{\varepsilon}_{\rho_t})\cdot\nabla\varphi_t  \,d\omega\\
&= -\int  \rho_t\,\nabla\varphi_t\cdot \nabla\bigl(\nabla(\omega\cdot\Omega^{\varepsilon}_{\rho_t})\cdot\nabla\varphi_t \bigr) \,d\omega\\
&=-\int\rho_t\,\nabla\varphi_t \,D^{2}(\omega\cdot\Omega^{\varepsilon}_{\rho_t})\,\nabla\varphi_t \,d\omega+\int\rho_t  \,\nabla(\omega\cdot \Omega^{\varepsilon}_{\rho_t})\cdot\nabla\Big(\frac{|\nabla\varphi_t|^2}{2}\Big)\,d\omega.
\end{aligned}
\end{align*}
Similarly we have
\begin{align*}
\begin{aligned}
I_{22}&= \int\rho_t\, \nabla\Big(\frac{|\nabla\varphi_t|^2}{2}\Big)\cdot\Omega^{\varepsilon}_{\rho_t} \,d\omega\\
&= \int\rho_t \, \nabla(\omega\cdot \Omega^{\varepsilon}_{\rho_t})\cdot\nabla\Big(\frac{|\nabla\varphi_t|^2}{2}\Big)\,d\omega,
\end{aligned}
\end{align*}
thus
\[
I_{21}+I_{22} = -\int\rho_t\,\nabla\varphi_t \,D^{2}(\omega\cdot\Omega^{\varepsilon}_{\rho_t})\,\nabla\varphi_t \,d\omega.
\]
Concerning $I_{23}$, since
\begin{align*}
\begin{aligned}
\partial_t\Omega^{\varepsilon}_{\rho_t}&= \frac{\partial_{t}J_{\rho_t}}{\sqrt{|J_{\rho_t}|^{2}+\varepsilon}}
- \frac{\Omega^{\varepsilon}_{\rho_t}}{\sqrt{|J_{\rho_t}|^{2}+\varepsilon}}\,\Omega^{\varepsilon}_{\rho_t}\cdot \partial_{t}J_{\rho_t}\\
&= \frac{1}{\sqrt{|J_{\rho_t}|^{2}+\varepsilon}} \Big( -\int \omega\,\nabla\cdot(\rho_t\nabla\varphi_t) \,d\omega + \Omega^{\varepsilon}_{\rho_t}\int\Omega^{\varepsilon}_{\rho_t}\cdot\omega\,\nabla\cdot(\rho_t\nabla\varphi_t)  \,d\omega  \Big)\\
&= \frac{1}{\sqrt{|J_{\rho_t}|^{2}+\varepsilon}} \Big( \int \rho_t\,\nabla\varphi_t \,d\omega -\Omega^{\varepsilon}_{\rho_t}\int \nabla(\Omega^{\varepsilon}_{\rho_t}\cdot\omega)\cdot\nabla\varphi_t \,\rho_t\,  d\omega  \Big)\\
&= \frac{1}{\sqrt{|J_{\rho_t}|^{2}+\varepsilon}} \Big( \int \rho_t\,\nabla\varphi_t \,d\omega -\Omega^{\varepsilon}_{\rho_t}\int \Omega^{\varepsilon}_{\rho_t}\cdot\nabla\varphi_t\, \rho_t\,  d\omega  \Big),\\
\end{aligned}
\end{align*}
we have
\begin{align*}
\begin{aligned}
I_{23}&=- \frac{1}{\sqrt{|J_{\rho_t}|^{2}+\varepsilon}}\int\rho_t \, \nabla\varphi\cdot\Big( \int \rho_t\,\nabla\varphi_t \,d\omega - \Omega^{\varepsilon}_{\rho_t}\int\Omega^{\varepsilon}_{\rho_t}\cdot\nabla\varphi_t \,\rho_t  \,d\omega  \Big) \,d\omega\\
&=- \frac{1}{\sqrt{|J_{\rho_t}|^{2}+\varepsilon}}\bigg|\int\nabla\varphi_t\,\rho_t \,d\omega\bigg|^{2}
+ \frac{1}{\sqrt{|J_{\rho_t}|^{2}+\varepsilon}}\bigg(\int\Omega^{\varepsilon}_{\rho_t}\cdot\nabla\varphi_t\,\rho_t \,d\omega\bigg)^{2}.
%&=- \frac{1}{\sqrt{|J_{\rho_t}|^{2}+\varepsilon}}\bigg|\int\nabla\varphi_t\,\rho_t \,d\omega\bigg|^{2}+\frac{1}{\sqrt{|J_{\rho_t}|^{2}+\varepsilon}}\bigg(\int\nabla\varphi_t\,\rho_t \,d\omega\bigg)\big(\Omega_{\rho_t}^{\varepsilon}\otimes\Omega_{\rho_t}^{\varepsilon}\big)\bigg(\int\nabla\varphi_t\,\rho_t \,d\omega\bigg).
\end{aligned}
\end{align*}
Recalling that $\rho_0=\rho$ and $\varphi_0=\varphi$, this completes the proof of \eqref{Hessian}.
\end{proof}

\subsection{Minimizing movements for the free energy}
To prove existence of solutions to the regularized problem,
we use the time-discrete scheme by Jordan, Kindeleherer and Otto \cite{JKO} (see also \cite{Figalli-Gigli}).
Hence, in all this section, $\eps>0$ is fixed and, to simplify the notation, we shall not explicitly show the dependence on it.

Given a time step $\tau>0$, for a given initial data $\rho_{0}\in\mathcal{P}(\bbs^{d-1})$ we set 
\[
\rho_{0}^{\tau}=\rho_{0},
\]
and then recursively define $\rho_{n}^{\tau}$ as a minimizer of
\begin{align}
\begin{aligned}\label{scheme-W}
\sigma\mapsto \frac{W_{2}^{2}(\sigma,\rho_{n-1}^{\tau})}{2\tau}+\mathcal{E^{\varepsilon}}(\rho).
\end{aligned}
\end{align} 
The existence of a minimizer to \eqref{scheme-W} is guaranteed as follows.

\begin{lemma}\label{lem-step} For a given $\tau>0$ and $\rho\in\mathcal{P}(\bbs^{d-1})$, there exists
a minimum $\rho_{\tau}\in\mathcal{P}(\bbs^{d-1})$ of
\[
\sigma\rightarrow\frac{W_{2}^{2}(\sigma,\rho)}{2\tau}+\mathcal{E^{\varepsilon}}(\sigma)\label{minmov}.
\]
Furthermore, the optimal transport map $T$ sending $\rho_\tau\,d\omega$ onto $\rho\,d\omega$ is given by
\beq\label{T-map}
T(\omega)=\exp_{\omega}\big[\tau\nabla_{\omega}\big(\log\rho_{\tau}-\omega\cdot\Omega^{\varepsilon}_{\rho_{\tau}}\big)\big].
\eeq
\end{lemma}
\begin{proof}
First of all, the existence of a minimum $\mu_{\tau}=\rho_{\tau}\,d\omega$ follows from the fact that 
$\mathcal{E^{\varepsilon}}$ is lower semicontinous and bounded from below thanks to 
\begin{align}
\begin{aligned}\label{E-lower}
\mathcal{E}^{\varepsilon}(\rho_{t})&\geq \min_{x>0}x\log x \int_{\bbs^{d-1}} \,d\omega -\sqrt{\Big|\int_{\bbs^{d-1}}\rho_t \,d\omega \Big|^2+\varepsilon}\\
%&\geq \min_{x>0}x\log x \int_{\bbs^{d-1}} \,d\omega -\sqrt{1+\varepsilon}\\
&\geq |\bbs^{d-1}|\, e^{-1}\log e^{-1} -\sqrt{1+\varepsilon}.
\end{aligned}
\end{align}
To show \eqref{T-map}, let $\varphi$ be a $d^2/2$-convex function such that $\exp_{\omega}(\nabla_{\omega}\varphi)$ is the optimal map 
sending $\rho_\tau\,d\omega$ onto $\rho\,d\omega.$
For any smooth vector field $X$ on $\mathbb{S}^{d-1}$, we set
\[
\mu_t:=\exp(tX)_{\#}\rho_{\tau} \,d\omega.
\]
Using the minimality of 
$\rho_{\tau}$, we get
\[
 \mathcal{E}^{\varepsilon}(\mu_t)-\mathcal{E}^{\varepsilon}(\rho_{\tau})+
 \frac{W_{2}^{2}(\mu_t, \rho)-
 W_{2}^{2}(\rho_{\tau}, \rho)}{2\tau}\geq0.
\]
Then we use Proposition \ref{prop-W} and \eqref{direct-E} to obtain
\begin{align*}
\begin{aligned}
&\int_{\mathbb{S}^{d-1}}\nabla_{\omega}\big(\log\rho_{\tau}-\omega\cdot\Omega^{\varepsilon}_{\rho_{\tau}}\big)\cdot X(\omega)\,\rho \,d\omega
 -\frac{1}{\tau}\int_{\mathbb{S}^{d-1}}\nabla_{\omega}\varphi\cdot
 X(\omega)\,\rho \,d\omega\\
&\quad\geq  \limsup_{t\rightarrow0}\Big( \mathcal{E}^{\varepsilon}(\mu_t)-\mathcal{E}^{\varepsilon}(\rho_{\tau})+
 \frac{W_{2}^{2}(\mu_t, \rho)-
 W_{2}^{2}(\rho_{\tau}, \rho)}{2\tau} \Big)\ge 0.
\end{aligned}
\end{align*}
Exchanging $X$ with $-X$, this yields
\[
\int_{\mathbb{S}^{d-1}}\tau\,\nabla_{\omega}\big(\log\rho_{\tau}-\omega\cdot\Omega^{\varepsilon}_{\rho_{\tau}}\big)\cdot X(\omega)\,\rho \,d\omega = \int_{\mathbb{S}^{d-1}}\nabla_{\omega}\varphi\cdot
 X(\omega)\,\rho \,d\omega,
\]
and since $X$ is arbitrary we get
\[
\nabla_{\omega}\varphi=\tau\,\nabla_{\omega}\big(\log\rho_{\tau}-\omega\cdot\Omega^{\varepsilon}_{\rho_{\tau}}\big),
\]
which proves \eqref{T-map}.
\end{proof}

\subsection{Existence of the regularized equation \eqref{eqrob}}
Using the sequence of minimizers defined in the previous section, we define the discrete solution $t\mapsto \rho^{\tau}(t)$ by
$$
\rho^{\tau}(t):=\rho_{n}^{\tau},\qquad \text{for } t\in[n\tau,(n+1)\tau).
$$
We show now the existence of weak solutions to \eqref{eqrob} as a limit of the discrete solutions $\rho^{\tau}$ as $\tau\to 0$.
\begin{proposition} \label{prop-exist} 
Assume $\rho_0\in\mathcal{P}(\bbs^{d-1})$ with $\int_{\bbs^{d-1}}\rho_0\log \rho_0 \,d\omega<\infty$. Then, for any sequence $\tau_{k}\downarrow0$,
up to a subsequence $\rho^{\tau_{k}}(t)$ converges to some limit $\rho(t)$ locally uniformly in time. The limit $t\mapsto\rho(t)$
belongs to $L_{loc}^{2}([0,\infty),W^{1,1}(\bbs^{d-1}))$ and is a weak solution of \eqref{eqrob}.
\end{proposition}
\begin{proof}
Throughout the proof, we will use the following inequality for the sequence of minimizers $\rho_n^{\tau}$ of \eqref{scheme-W},
\begin{equation}
 \label{e-ineq}
\frac{1}{2}\sum_{i=n}^{m-1}\frac{W_{2}^{2}(\rho_{i+1}^{\tau},\rho_{i}^{\tau})}{\tau}+\frac{\tau}{2}\sum_{i=n}^{m-1}|\nabla \mathcal{E}^{\varepsilon}(\rho_{i}^{\tau})|^{2}\leq \mathcal{E}^{\varepsilon}(\rho_{m}^{\tau})-\mathcal{E}^{\varepsilon}(\rho_{n}^{\tau}),\qquad \mbox{for any}~ n< m,
\end{equation}
referring to \cite[Lemma 3.2.2]{Ambrosio-Gigli-Savare} for its proof.\\ 
Since $\mathcal{E}^{\varepsilon}(\rho_{m}^{\tau})\le \mathcal{E}^{\varepsilon}(\rho_0)$ for all $m$ and
$\mathcal{E}^{\varepsilon}(\rho_{n}^{\tau})$ bounded from below due to \eqref{E-lower}, we have
\beq\label{e-1}
 \mathcal{E}^{\varepsilon}(\rho_{m}^{\tau})-\mathcal{E}^{\varepsilon}(\rho_{n}^{\tau})\leq {\mathcal{E}^{\varepsilon}(\rho_{0})}+\sqrt{1+\eps}.
\eeq
Let $\{\tau_{k}\}_{k\in\mathbb{N}}$ be a sequence converging to $0$. Then, for any $n< m$,
 \beq\label{useful}
\frac{1}{2}\sum_{i=n}^{m-1}\frac{W_{2}^{2}(\rho_{i+1}^{\tau_{k}},\rho_{i}^{\tau_{k}})}{\tau_k}\leq {\mathcal{E}^{\varepsilon}(\rho_{0})}+\sqrt{2}
 \eeq
 (recall that $\eps \leq 1$).
 
Notice that 
 \begin{equation}
 \label{eq:Eeps}
\mathcal{E}^{\varepsilon}(\rho_{0}) \le \int_{\bbs^{d-1}}\rho_0\log \rho_0\,d\omega<\infty. 
 \end{equation}
Also, it follows by Jensen's inequality that
\[
\frac{W_{2}^{2}(\rho_{m}^{\tau},\rho_{n}^{\tau})}{(m-n)^2}\le \Big(\frac{\sum_{i=n}^{m-1}W_{2}(\rho_{i+1}^{\tau},\rho_{i}^{\tau})}{m-n}\Big)^2 \le \frac{\sum_{i=n}^{m-1}W_{2}^2(\rho_{i+1}^{\tau},\rho_{i}^{\tau})}{m-n}\le 2\tau ({\mathcal{E}^{\varepsilon}(\rho_{0})}+\sqrt{2}).
\]
Hence, setting $n=[\frac{s}{\tau}]$ and $m=[\frac{t}{\tau}]$ for any $0\le s<t$, we have
 \begin{equation}
 \label{d-equicont}
  W_{2}(\rho^{\tau_{k}}(t),\rho^{\tau_{k}}(s))\leq\sqrt{2(\mathcal{E}^{\varepsilon}(\rho_{0})+\sqrt{2})\,[t-s+\tau_{k}]}.
 \end{equation}
This equicontinuity estimate combined the compactness of $(\mathcal{P}(\bbs^{d-1}),W_{2})$ implies that, up to a subsequence, 
\beq\label{converge}
\rho^{\tau_{k}}(t)~\mbox{converges to some limit}~\rho(t)~\mbox{in} ~(\mathcal{P}(\mathbb{S}^{d-1}),W_{2})~\mbox{ locally uniformly in} ~t\ge0. 
\eeq
We now show that $t\mapsto\rho(t)$ is a weak solution of 
 \eqref{eqrob}. For $n\in\mathbb{N}$, by \eqref{scheme-W} and Lemma \ref{lem-step}, we have
\[
\Big(\exp_{\omega}\big(\tau_{k}\nabla_{\omega}\big(\log\rho^{\tau_{k}}_{n+1}-\omega\cdot\Omega^{\varepsilon}_{\rho^{\tau_{k}}_{n+1}}\big)\big)\Big)_{\#}\rho^{\tau_{k}}_{n+1}\,d\omega = \rho^{\tau_{k}}_{n}\,d\omega.
\]
Thus, for any $\varphi\in C^{\infty}(\mathbb{S}^{d-1})$,
\[
\int_{\bbs^{d-1}}\varphi(\omega)\,(\rho_{n+1}^{\tau_{k}}-\rho_{n}^{\tau_{k}})\,d\omega
=\int_{\bbs^{d-1}}\Big(\varphi(\omega)
-\varphi \big(\exp_{\omega}(\tau_{k}\nabla_{\omega}(\log\rho^{\tau_{k}}_{n+1}-\omega\cdot\Omega^{\varepsilon}_{\rho^{\tau_{k}}_{n+1}}))\big)\Big)
\,\rho^{\tau_{k}}_{n+1} \,d\omega.
\]
Using, for each $\omega \in \bbs^{d-1}$, the Taylor formula along the geodesic $s\mapsto\exp_{\omega}(s\tau_{k}\nabla_{\omega}(\log\rho^{\tau_{k}}_{n+1}-\omega\cdot\Omega^{\varepsilon}_{\rho^{\tau_{k}}_{n+1}}))$, we have
\begin{align}
\begin{aligned}\label{sum-eq}
\int_{\bbs^{d-1}}\varphi(\omega)\,(\rho_{n+1}^{\tau_{k}}-\rho_{n}^{\tau})\,d\omega
&=-\int_{\bbs^{d-1}}\tau_{k}\,\nabla_{\omega}\varphi(\omega)\cdot
\nabla_{\omega}\big[\log\rho_{n+1}^{\tau_{k}}
-\omega\cdot\Omega^{\varepsilon}_{\rho^{\tau_{k}}_{n+1}}\big]\,
\rho_{n+1}^{\tau_{k}} \,d\omega
+R(n,\tau_{k})\\
&=\int_{\bbs^{d-1}}\tau_{k}\big[\Delta_{\omega}\varphi
+\nabla_{\omega}\varphi\cdot \nabla_{\omega}(\omega\cdot\Omega^{\varepsilon}_{\rho^{\tau_{k}}_{n+1}})\big]
\rho_{n+1}^{\tau_{k}} \,d\omega+R(n,\tau_{k}),
\end{aligned}
\end{align}
where the remainder term $R(n,\tau_k)$ can be estimated by
\begin{align}
\begin{aligned}\label{remainder}
R(n,\tau_k)&\leq\|D_{\omega}^{2}\varphi\|_{L^{\infty}(\bbs^{d-1})}\int_{S^{d-1}} d^{2}\Big(\omega,\exp_{\omega}(\tau_{k}\nabla_{\omega}\big[\log\rho_{n+1}^{\tau_{k}}-\omega\cdot\Omega_{\rho^{\tau_{k}}_{n+1}}\big]\Big)\,\rho_{n+1}^{\tau_{k}}\,
d\omega\\
&=\|D_{\omega}^{2}\varphi\|_{L^{\infty}(\bbs^{d-1})}\,
W_{2}^{2}(\rho_{n+1}^{\tau_{k}},\rho_{n}^{\tau_k}),
\end{aligned}
\end{align}
For any $0\leq t<s $, we sum up \eqref{sum-eq} from $l:=[\frac{t}{\tau_{k}}]$ to $m:=[\frac{s}{\tau_{k}}]$ to get
\begin{align*}
\begin{aligned}
 \int_{\bbs^{d-1}}\varphi\,(\rho^{\tau_{k}}(s)-\rho^{\tau_{k}}(t))\,d\omega
 &\le \sum_{n=l}^{m}\int_{\bbs^{d-1}}\tau_{k}\big[\Delta_{\omega}\varphi
+\nabla_{\omega}\varphi\cdot \nabla_{\omega}(\omega\cdot\Omega^{\varepsilon}_{\rho^{\tau_{k}}_{n+1}})\big]\,
\rho_{n+1}^{\tau_{k}} \,d\omega \\
&\quad +\sum_{n=l}^{m}R(n,\tau_k)\\
 &=\int_{(l+1)\tau_k}^{(m+2)\tau_k}\int_{\bbs^{d-1}}\big[\Delta_{\omega}\varphi
+\nabla_{\omega}\varphi\cdot \nabla_{\omega}(\omega\cdot\Omega^{\varepsilon}_{\rho^{\tau_{k}}(r)})\big]\,
\rho^{\tau_{k}}(r) \,d\omega \,dr \\
&\quad +\sum_{n=l}^{m}R(n,\tau_k).
\end{aligned}
\end{align*}
Letting $\tau_{k}\to0$, \eqref{converge} implies
\[
 \int_{\bbs^{d-1}}\varphi\,(\rho^{\tau_{k}}(s)-\rho^{\tau_{k}}(t))\,d\omega \to \int_{\bbs^{d-1}}\varphi\,(\rho(s)-\rho(t))\,d\omega.
\]
Since $J_{\rho^{\tau_{k}}}\to J_{\rho}$, we have 
\begin{align*}
\begin{aligned}
|\Omega^{\varepsilon}_{\rho^{\tau_{k}}}-\Omega^{\eps}_{\rho}| &\le \frac{\Big| J_{\rho^{\tau_{k}}}(\sqrt{|J_{\rho}|^2+\eps}- \sqrt{|J_{\rho^{\tau_{k}}}|^2+\eps})+ \sqrt{|J_{\rho^{\tau_{k}}}|^2+\eps}(J_{\rho^{\tau_{k}}}-J_{\rho}) \Big|}
{\sqrt{|J_{\rho^{\tau_{k}}}|^2+\eps}\sqrt{|J_{\rho}|^2+\eps}}\\
& \le \frac{1}{\eps}\Big(|J_{\rho^{\tau_{k}}}-J_{\rho}|+\Big|\sqrt{|J_{\rho^{\tau_{k}}}|^2+\eps}- \sqrt{|J_{\rho}|^2+\eps}\Big| \Big)\\
&\to 0,
\end{aligned}
\end{align*}
which implies that, for all $r$,
\begin{align*}
\begin{aligned}
&\int_{\bbs^{d-1}}\big[\Delta_{\omega}\varphi
+\nabla_{\omega}\varphi\cdot\nabla_{\omega}(\omega\cdot\Omega^{\varepsilon}_{\rho^{\tau_{k}}(r)})\big]\,
\rho^{\tau_{k}}(r) \,d\omega\\
&\quad\to  \int_{\bbs^{d-1}}\big[\Delta_{\omega}\varphi
+\nabla_{\omega}\varphi\cdot\nabla_{\omega}(\omega\cdot\Omega^{\varepsilon}_{\rho(r)})\big]\,
\rho(r) \,d\omega.
\end{aligned}
\end{align*}
Moreover since $\int_{\bbs^{d-1}}\big[\Delta_{\omega}\varphi
+\nabla_{\omega}\varphi\cdot \nabla_{\omega}(\omega\cdot\Omega^{\varepsilon}_{\rho^{\tau_{k}}(r)})\big]\,
\rho^{\tau_{k}}(r) \,d\omega$ is uniformly bounded, the dominated convergence theorem yields
\begin{align*}
\begin{aligned}
&\int_{(l+1)\tau_k}^{(m+2)\tau_k}\int_{\bbs^{d-1}}\big[\Delta_{\omega}\varphi
+\nabla_{\omega}\varphi\cdot\nabla_{\omega}(\omega\cdot\Omega^{\varepsilon}_{\rho^{\tau_{k}}(r)})\big]\,
\rho^{\tau_{k}}(r) \,d\omega \,dr\\
&\quad\to \int_{t}^s \int_{\bbs^{d-1}}\big[\Delta_{\omega}\varphi
+\nabla_{\omega}\varphi\cdot \nabla_{\omega}(\omega\cdot\Omega^{\varepsilon}_{\rho(r)})\big]\,
\rho(r) \,d\omega \,dr.
\end{aligned}
\end{align*}
On the other hand, since \eqref{remainder} and \eqref{useful} give
\begin{align*}
\begin{aligned}
\sum_{n=l}^{m}R(n,\tau_k) &\leq C\,\sum_{n=l}^{m} W_{2}^{2}(\rho_{n+1}^{\tau_{k}},\rho_{n}^{\tau_k})\\
&\leq C(\mathcal{E}^{\eps}(\rho_{0})+\sqrt2)\tau_{k} \to 0,
\end{aligned}
\end{align*}
we have shown that $0\leq t<s $,
\[
\int_{\bbs^{d-1}}\varphi\,(\rho(s)-\rho(t))\,d\omega =\int_{t}^s \int_{\bbs^{d-1}}\big[\Delta_{\omega}\varphi
+\nabla_{\omega}\varphi\cdot \nabla_{\omega}(\omega\cdot\Omega^{\varepsilon}_{\rho(r)})\big]\,
\rho(r) \,d\omega \,dr,
\]
which provides the weak formulation of $\eqref{eqrob}$. Moreover, thanks to \eqref{d-equicont} and \eqref{eq:Eeps},
 \begin{equation}
 \label{eq:equicont}
W_2(\rho(t),\rho(s)) \leq \sqrt{2 \int_{\bbs^{d-1}}\rho_0\log \rho_0\,d\omega}\,\sqrt{t-s},
\end{equation}
hence
$t\mapsto \rho(t)$ is weakly continuous and $\rho$ is a weak solution to $\eqref{eqrob}$.

It remains to show that $\rho\in L_{loc}^{2}([0,\infty),W^{1,1}(\bbs^{d-1}))$. Using again \eqref{e-ineq} and \eqref{e-1} we see that, for any $0\le t<s$,
\[
\int_t^s |\nabla\mathcal{E}^{\varepsilon}(\rho^{\tau_{k}}(t))|^{2}dt\leq 2(\mathcal{E}^{\eps}(\rho_{0})+\sqrt2),
\]
which together with \eqref{E-slope} yields
\begin{align}
\begin{aligned}\label{e-com}
2(\mathcal{E}^{\eps}(\rho_{0})+\sqrt2)&\ge \int_t^s\int_{\bbs^{d-1}}|\nabla_{\omega}(\log\rho^{\tau_{k}}-\omega\cdot\Omega^{\varepsilon}_{\rho^{\tau_{k}}})|^{2}\rho^{\tau_{k}} ~d\omega \,dt \\
 &\ge \int_t^s\int_{\bbs^{d-1}}|\nabla_{\omega}\log\rho^{\tau_{k}}|^2\rho^{\tau_{k}} \,d\omega \,dt
 -2\int_t^s\int_{\bbs^{d-1}}\nabla_{\omega}\log\rho^{\tau_{k}}\cdot\nabla_{\omega}(\omega\cdot\Omega^{\varepsilon}_{\rho^{\tau_{k}}})\,\rho^{\tau_{k}} \,d\omega \,dt\\
&\ge \frac{1}{2}\int_t^s\int_{\bbs^{d-1}}|\nabla_{\omega}\log\rho^{\tau_{k}}|^2\rho^{\tau_{k}} \,d\omega \,dt
 -2\int_t^s\int_{\bbs^{d-1}}|\nabla_{\omega}(\omega\cdot\Omega^{\varepsilon}_{\rho^{\tau_{k}}})|^2\rho^{\tau_{k}} \,d\omega \,dt.
\end{aligned}
\end{align}
Since $|\nabla_{\omega}(\omega\cdot\Omega^{\varepsilon})|=|\bbp_{\omega^{\perp}}\Omega|\le 1$, we have
\begin{align*}
\int_t^s\int_{\bbs^{d-1}}|\nabla_{\omega}\sqrt{\rho^{\tau_{k}}}|^2 d\omega \,dt
&=\frac{1}{2}\int_t^s\int_{\bbs^{d-1}}|\nabla_{\omega}\log\rho^{\tau_{k}}|^2\rho^{\tau_{k}} \,d\omega \,dt\\
&\le 2(\mathcal{E}^{\eps}(\rho_{0})+\sqrt2) + 2(s-t),
\end{align*}
which implies that $\sqrt{\rho^{\tau_{k}}}$ is uniformly bounded in $L^{2}_{loc}([0,+\infty),H^{1}(\bbs^{d-1}))$.\\ 
Therefore, letting $\tau_k \to 0$, we get
\beq\label{sqrt}
\sqrt{\rho}\in L^{2}_{loc}([0,+\infty),H^{1}(\bbs^{d-1})),
\eeq 
that combined with H\"older inequality implies that 
${\rho}\in L^{2}_{loc}([0,+\infty),W^{1,1}(\bbs^{d-1}))$. 
\end{proof}

\subsection{Uniqueness}
The following results provide the stability estimates for weak solutions to \eqref{eqrob}, thus their uniqueness. We shall revisit the arguments of the proof to show the stability and uniqueness of weak solutions to \eqref{main} in Section 5. 
\begin{proposition}\label{prop-unique} (Uniqueness and stability). 
Assume $\rho_0,\bar\rho_0\in\mathcal{P}(\bbs^{d-1})$ satisfy \eqref{ini-assume}.  Let $\rho^{\varepsilon}, \bar{\rho}^{\varepsilon}$ be solutions of \eqref{eqrob} with corresponding initial datas $\rho_0,\bar\rho_0$. Then for all $t>0$,
\beq\label{eps-uni}
W_{2}(\rho^{\varepsilon}(t),\bar{\rho}^{\varepsilon}(t))\leq e^{\lambda t}W_{2}(\rho_{0},\bar{\rho}_{0}),
\eeq
where $\lambda:=(1+\varepsilon^{-1/2})-(d-2)$.
\end{proposition}

\begin{proof}
For a fixed time $t>0$, let $\varphi_0$ be a $d^2/2$-convex function such that $\exp_{\omega}(\nabla_{\omega}\varphi_0)$ is the optimal map 
sending $\rho^\eps(t)\,d\omega$ onto $\bar\rho^\eps(t)\,d\omega$, and 
consider the curve $[0,1]\ni r\mapsto\alpha_{r}\,d\omega$ of absolutely continuous measures defined by
\[
\alpha_{r} \,d\omega=\exp_{\omega}(r\nabla_{\omega}\varphi_{0})_{\#} \rho^{\varepsilon}(t)\,d\omega
\]
(the absolute continuity of $\alpha_r$ follows, for instance, from \cite[Section 5]{F-F}).
Then the curve $r\mapsto\alpha_{r} \,d\omega$ is the unique geodesic in $(\mathcal{P}(\bbs^{d-1}),W_{2})$ connecting $\alpha_{0}=\rho^{\varepsilon}(t)$ to $\alpha_{1}=\bar\rho^{\varepsilon}(t)$ (see for example \cite[Corollary 3.22]{user}).\\
For each $r\in[0,1]$, let $\varphi_r$ be a $d^2/2$-convex function such that $\exp_{\omega}(\nabla_{\omega}\varphi_r)$ is the optimal map 
sending $\alpha_{r} \,d\omega$ onto $\bar\rho^{\varepsilon}(t)\,d\omega.$
Similarly, the curve $s\mapsto\alpha_{r,s} \,d\omega$ defined by
 \beq\label{exp-alpha}
\alpha_{r,s} \,d\omega = \exp_{\omega}(s\nabla\varphi_{r})_{\#}\alpha_{r} \,d\omega,
\eeq
and it is the unique geodesic in $(\mathcal{P}(\bbs^{d-1}),W_{2})$ connecting $\alpha_{r,0}=\alpha_{r}$ to $\alpha_{r,1}=\bar\rho^{\varepsilon}(t)$. Notice that it follows from the uniqueness of the geodesics that, for all $r,s\in[0,1]$,
\[
 \alpha_{r+(1-r)s}=\alpha_{r,s}.
\]
Now, applying \eqref{Hessian} in Lemma \ref{lem-E} to \eqref{exp-alpha}, we estimate the second derivative of $\mathcal{E^{\varepsilon}}$ by Wasserstein distance as
\begin{align}
\begin{aligned}\label{second-1}
\frac{d^{2}}{dh^{2}}\bigg|_{h=r}\mathcal{E^{\varepsilon}}(\alpha_{h})&=\frac{d^{2}}{dh^{2}}\bigg|_{h=0}\mathcal{E^{\varepsilon}}(\alpha_{r,\frac{h}{1-r}})\\
&=\frac{1}{(1-r)^2}\frac{d^{2}}{ds^{2}}\bigg|_{s=0}\mathcal{E^{\varepsilon}}(\alpha_{r,s})\\
&\geq-\frac{\lambda}{(1-r)^{2}}\int_{\bbs^{d-1}}|\nabla\varphi_{r}|^{2}\alpha_{r} \,d\omega\\
&=-\lambda\,\frac{W_{2}^{2}(\alpha_{r},\bar\rho^{\varepsilon}(t))}{(1-r)^{2}}\\
&=-\lambda\, W_{2}^{2}(\rho^{\varepsilon}(t),\bar\rho^{\varepsilon}(t)),
\end{aligned}
\end{align}
where $\lambda:=(1+\varepsilon^{-1/2})-(d-2)$ (see \eqref{see}).\\
Since, by Taylor formula along the geodesic $r\mapsto\alpha_{r} \,d\omega$, 
\[
 \mathcal{E^{\varepsilon}}(\alpha_{1})=\mathcal{E^{\varepsilon}}(\alpha_{0})+\frac{d}{dr}\bigg|_{r=0}\mathcal{E^{\varepsilon}}(\alpha_{r})+\int_{0}^{1}(1-r)\frac{d^{2}}{dr^{2}}\mathcal{E^{\varepsilon}}(\alpha_{r})\,dr,
\]
we use \eqref{direct-E} and \eqref{second-1} to have
\[
 \mathcal{E^{\varepsilon}}(\bar\rho^{\varepsilon}(t))\geq\mathcal{E^{\varepsilon}}(\rho^{\varepsilon}(t))+\int_{\bbs^{d-1}}\nabla\varphi_{0}\cdot\nabla(\log\rho^{\varepsilon}(t)-\omega\cdot\Omega_{\rho^{\varepsilon}(t)})\,\rho^{\varepsilon}(t) \,d\omega-\frac{\lambda}{2}\,W_{2}^{2}(\rho^{\varepsilon}(t),\bar\rho^{\varepsilon}(t)).
\]
Similarly, applying the above arguments to the $d^2/2$-convex function $\overline{\varphi}_{0}$ satisfying
\[
 \rho^{\varepsilon}(t)\,d\omega=\exp_{\omega}(\nabla\overline{\varphi}_{0})_{\#}\overline{\rho}^{\varepsilon}(t)\,d\omega,
\]
we have
\[
 \mathcal{E^{\varepsilon}}(\rho^{\varepsilon}(t))\geq\mathcal{E^{\varepsilon}}(\bar\rho^{\varepsilon}(t))+\int_{\bbs^{d-1}}\nabla\bar\varphi_{0}\cdot
 \nabla(\log\bar\rho^{\varepsilon}(t)-\omega\cdot\Omega_{\bar\rho^{\varepsilon}(t)})\,\bar\rho^{\varepsilon}(t) \,d\omega-\frac{\lambda}{2}\,W_{2}^{2}(\rho^{\varepsilon}(t),\bar\rho^{\varepsilon}(t)),
\]
therefore 
\begin{multline}
\label{eq:lambda}
\int_{\bbs^{d-1}}\nabla\varphi_{0}\cdot\nabla(\log\rho^{\varepsilon}(t)-\omega\cdot\Omega_{\rho^{\varepsilon}(t)})\,\rho^{\varepsilon}(t) \,d\omega\\
+ \int_{\bbs^{d-1}}\nabla\bar\varphi_{0}\cdot\nabla(\log\bar\rho^{\varepsilon}(t)-\omega\cdot\Omega_{\bar\rho^{\varepsilon}(t)})\,\bar\rho^{\varepsilon}(t) \,d\omega
\le \lambda \,W_{2}^{2}(\rho^{\varepsilon}(t),\bar\rho^{\varepsilon}(t)).
\end{multline}
We now claim that
\begin{equation}
\label{eq:claim}
\begin{split}
\frac{d}{dt}W_{2}^{2}(\rho^{\varepsilon}(t),\bar\rho^{\varepsilon}(t))&= \int_{\bbs^{d-1}}\nabla\varphi_{0}\cdot\nabla(\log\rho^{\varepsilon}(t)-\omega\cdot\Omega_{\rho^{\varepsilon}(t)})\,\rho^{\varepsilon}(t) \,d\omega\\
&\quad+ \int_{\bbs^{d-1}}\nabla\bar\varphi_{0}\cdot\nabla(\log\bar\rho^{\varepsilon}(t)-\omega\cdot\Omega_{\bar\rho^{\varepsilon}(t)})
\,\bar\rho^{\varepsilon}(t) \,d\omega.
\end{split}
\end{equation}
Indeed, $\rho^{\eps}$ and $\bar\rho^{\eps}$ solve the continuity equation
\[
\partial_{t}\rho+\nabla_{\omega}\cdot(v[\rho]\rho)=0,
\]
where $v[\rho]:=\nabla_{\omega}(\omega\cdot\Omega^{\eps}_\rho-\log\rho)$ is a locally Lipschitz vector field. Moreover it follows
from \eqref{sqrt} that, for all $t<s$,
\[
\int_t^s \int_{\bbs^{d-1}} |v[\rho]|\,\rho \,d\omega\le C\,(s-t)\Big(1+\|\nabla\sqrt{\rho}\|_{L^2(\bbs^{d-1})}\Big)<\infty.
\]
Hence the hypotheses of \cite[Theorem 23.9]{Villani} are satisfied implying \eqref{eq:claim},
and combining it with 
\eqref{eq:lambda} yields
\[
\frac{d}{dt}W_{2}^{2}(\rho^{\varepsilon}(t),\bar\rho^{\varepsilon}(t)\le \lambda\,W_{2}^{2}(\rho^{\varepsilon}(t),\bar\rho^{\varepsilon}(t)),
\]
which completes the proof.
\end{proof}

\subsection{Properties of solutions to \eqref{eqrob}}
In the next lemma, we show that the momentum does not vanish for any finite time $t\in(0,\infty)$,
with an estimate independent of $\eps$.
\begin{lemma} \label{lem-moment}  
Let $\rho^{\varepsilon}$ be a solution of \eqref{eqrob}. Then, for all $t>0$,
\[
| J_{{\rho}^{\varepsilon}(t)} |^{2}\geq | J_{\rho_{0}}|^{2}e^{-2(d-1)t}.
\]
\end{lemma}
\begin{proof}
It follows from \eqref{eqrob} that
\begin{align*}
\begin{aligned} 
\frac{d}{dt}\frac{1}{2}| J_{{\rho}^{\varepsilon}}|^{2}&=J_{{\rho}^{\varepsilon}}\cdot\partial_{t}J_{{\rho}^{\varepsilon}}\\
&=J_{{\rho}^{\varepsilon}}\cdot \bigg(\int\omega\,\Delta\rho^{\varepsilon} \,d\omega-\int\omega\, \nabla\cdot(\rho^{\eps}\nabla(\omega\cdot\Omega_{{\rho}^{\varepsilon}}))\,d\omega\bigg)\\
&=J_{{\rho}^{\varepsilon}}\cdot \int\omega\,\Delta\rho^{\varepsilon} \,d\omega-\int J_{{\rho}^{\varepsilon}}\cdot \omega\,
\nabla\cdot(\rho^{\eps}\nabla(\omega\cdot\Omega_{{\rho}^{\varepsilon}}))\,d\omega \\
&=: I_1 + I_2.
\end{aligned}
\end{align*}
We use \eqref{form-0} to get
\begin{align*}
\begin{aligned} 
I_1 &= -J_{{\rho}^{\varepsilon}}\cdot \int\nabla\rho^{\varepsilon} \,d\omega\\
&= -(d-1)J_{{\rho}^{\varepsilon}}\cdot \int\omega \,\rho^{\varepsilon} \,d\omega\\
&=-(d-1)|J_{{\rho}^{\varepsilon}}|^2.
\end{aligned}
\end{align*}
Also, using \eqref{formula-0}, we have
\begin{align*}
\begin{aligned} 
I_2 &= \int \nabla(J_{{\rho}^{\varepsilon}}\cdot \omega) \cdot\nabla(\omega\cdot\Omega_{{\rho}^{\varepsilon}})\,\rho^{\eps}
d\omega\\
&= \int \frac{\rho^{\eps}}{\sqrt{|J_{{\rho}^{\varepsilon}}|^{2}+\varepsilon}}\, |\nabla(\omega\cdot J_{{\rho}^{\varepsilon}})|^2\, d\omega\\
&\ge 0.
\end{aligned}
\end{align*}
Thus 
\[
\frac{d}{dt}| J_{{\rho}^{\varepsilon}}|^{2} \le -2(d-1)|J_{{\rho}^{\varepsilon}}|^2,
\]
which completes the proof.
\end{proof}

\subsection{Proof of the existence in Theorem \ref{thm-exist}}
Let $\{\eps_{k}\}_{k\in\mathbb{N}}$ be a sequence converging to $0$. 
As a consequence of \eqref{eq:equicont} it follows 
that the sequence $\{\rho^{\eps_k}\}_{k\in\mathbb{N}}$ is equicontinous, so the compactness of $(\mathcal{P}(\bbs^{d-1}),W_{2})$ imply that up to a subsequence, $\rho^{\eps_{k}}(t)$ converges to some limit $\rho(t)$ in $(\mathcal{P}(\mathbb{S}^{d-1}),W_{2})$ uniformly in $t\ge0$. Then since $J_{\rho^{\eps_k}}\to J_{\rho}$, it follows from Lemma \ref{lem-moment} that for all $t>0$,
\beq\label{J-0}
| J_{{\rho}(t)} |^{2}\geq | J_{\rho_{0}}|^{2}e^{-2(d-1)t}.
\eeq
Therefore by the same arguments as the proof of Proposition \ref{prop-exist}, the limit $\rho$ is a weak solution to \eqref{main}. Moreover since a straightforward computation yields
\begin{align*}
\frac{d}{dt}\mathcal{E}^{0}(\rho)& =-\int_{\bbs^{d-1}}|\nabla_{\omega}(\log\rho-\omega\cdot\Omega_{\rho})|^2 \rho\,d\omega,
\end{align*}
the analogue of \eqref{e-com} with $\eps=0$ combined with \eqref{eq:Eeps} provide
\[
{\rho}\in L^{2}_{loc}([0,+\infty),W^{1,1}(\bbs^{d-1})).
\]
\qed

\section{Convergence towards equilibrium}
\setcounter{equation}{0}
In this section, we prove Theorem \ref{thm-converge}. We start with the following estimates on the difference between $\rho^{\varepsilon}$ and $M_{\Omega^{\eps}_{\rho^{\varepsilon}}}$.
\begin{lemma}\label{lem-equili} 
Let $C_M$ be as in \eqref{eq:CM}, and let $\rho^{\varepsilon}$ be a solution of \eqref{eqrob}
starting from $\rho_0$. Then, for all $t>0$,
\[
\|\rho^{\varepsilon}(t)-M_{\Omega^{\eps}_{\rho^{\varepsilon}(t)}}\|_{L^{1}(\mathbb{S}^{d-1})}\leq e^{-C_1t}\Big(\int_{\bbs^{d-1}}\rho_0\log\rho_0\,d\omega +1-\log C_M\Big)+\sqrt{\eps}.
\]
where
\[
C_1:=\frac{2(d-2)}{e^{2}}.%, \quad C_2:= \frac{\varepsilon (d-1)}{C_1(|J_{\rho_{0}}|^{2}e^{-2(d-2)t}+\varepsilon)^{3/2}} (1-e^{-C_1t}).
\]
\end{lemma}
\begin{proof}
First of all, for each measure $\rho^{\eps}$, we denote its relative entropy with respect
to the probability measure $M_{\Omega^{\eps}_{\rho^{\varepsilon}}}(\omega)\,d\omega =C_M e^{\omega\cdot\Omega^{\eps}_{\rho^{\varepsilon}}}d\omega$ by 
\begin{align*}
\begin{aligned} 
H(\rho^{\eps}\mid M_{\Omega^{\eps}_{\rho^{\varepsilon}}} )=\int_{\bbs^{d-1}}\rho^{\eps}\log\Big(\frac{\rho^{\eps}}{M_{\Omega^{\eps}_{\rho^{\varepsilon}}} }\Big)\,d\omega,
\end{aligned}
\end{align*}
which can also be rewritten as
\[
H(\rho^{\eps}\mid M_{\Omega^{\eps}_{\rho^{\varepsilon}}} )=\int_{\bbs^{d-1}}\rho^{\eps}\log\rho^{\eps} \,d\omega-\int_{\bbs^{d-1}}
\omega\cdot\Omega^{\eps}_{{\rho}^{\eps}}\,\rho^{\eps} \,d\omega- \log C_M.
\]
Since
\begin{align*}
\begin{aligned} 
\int_{\bbs^{d-1}}\omega\cdot\Omega^{\eps}_{{\rho}^{\eps}}\,\rho^{\eps} \,d\omega &=
\frac{J_{\rho^{\eps}}}{\sqrt{|J_{\rho^{\eps}}|^{2}+\varepsilon}}\cdot\int\omega\,\rho^{\eps} \,d\omega=\frac{|J_{\rho^{\eps}}|^{2}}{\sqrt{|J_{\rho^{\eps}}|^{2}+\varepsilon}}\\
&=\sqrt{|J_{\rho^{\eps}}|^{2}+\varepsilon}-\frac{\varepsilon}{\sqrt{|J_{\rho^{\eps}}|^{2}+\varepsilon}},
\end{aligned}
\end{align*}
we have
\beq\label{relation-1}
\mathcal{E}^{\varepsilon}(\rho^{\eps})=H(\rho^{\eps}\mid M_{\Omega^{\varepsilon}_{\rho^{\eps}}})-\frac{\varepsilon}{\sqrt{|J_{\rho^{\eps}}|^{2}+\varepsilon}} + \log C_M.
\eeq
We now set $\alpha:=|\bbs^{d-1}|^{-1}$, and 
regard the measure $M_{\Omega^{\eps}_{\rho^{\varepsilon}}} \,d\omega$ as a bounded
perturbation of the constant probability measure $\alpha \,d\omega$, i.e.,
\[
M_{\Omega^{\eps}_{\rho^{\varepsilon}}} = e^{\omega\cdot\Omega^{\eps}_{\rho^{\varepsilon}}- \log C_M}=e^{\omega\cdot\Omega^{\eps}_{\rho^{\varepsilon}}- \log C_M-\log\alpha}\alpha,
\]
where 
\beq\label{osc}
{\rm osc}(\omega\cdot\Omega^{\eps}_{\rho^{\varepsilon}}- \log C_M-\log\alpha)\leq2.
\eeq
Since the Ricci curvature tensor of $\bbs^{d-1}$ is $(d-2)I_{d}$ and $d \geq 3$, the logarithmic Sobolev inequality \cite{B-E} implies
\[
H(\rho^{\eps}\mid \alpha)\leq\frac{1}{2(d-2)}\int\bigg|\nabla\log\frac{\rho^{\eps}}{\alpha}\bigg|^{2}\rho^{\varepsilon} \,d\omega.
\]
Thus,
since the logarithmic Sobolev inequality  is stable under bounded perturbations (see for instance \cite{H-S, Otto-Villani}),
  it follows from \eqref{osc} that
\beq\label{ineq-1}
H(\rho^{\eps}\mid M_{\Omega^{\eps}_{\rho^{\varepsilon}}})\leq\frac{e^{2}}{2(d-2)}\int\bigg|\nabla\log\frac{{\rho}^{\varepsilon}}{M_{\Omega^{\eps}_{\rho^{\varepsilon}}}}\bigg|^{2}\rho^{\varepsilon} \,d\omega.
\eeq
Therefore, since \eqref{relation-1} yields
\begin{align*}
\begin{aligned} 
\frac{d}{dt}\mathcal{E}^{\varepsilon}(\rho^{\eps})=\frac{d}{dt}H(\rho^{\eps}\mid M_{\Omega^{\varepsilon}_{\rho^{\eps}}})-\frac{d}{dt}\frac{\varepsilon}{\sqrt{|J_{\rho^{\eps}}|^{2}+\varepsilon}},
\end{aligned}
\end{align*} 
and we see
\begin{align}
\begin{aligned}\label{ineq-2}
\frac{d}{dt}\mathcal{E}^{\varepsilon}({\rho}^{\varepsilon})& =-\int |\nabla(\log{\rho}^{\varepsilon}- \omega\cdot\Omega^{\eps}_{\rho^{\varepsilon}})|^2 {\rho}^{\varepsilon} \,d\omega\\
&=-\int\bigg|\nabla\log\frac{{\rho}^{\varepsilon}}{M_{\Omega^{\eps}_{\rho^{\varepsilon}}}}\bigg|^{2}\rho^{\varepsilon} \,d\omega,
\end{aligned}
\end{align}
it follows from \eqref{ineq-1} and \eqref{ineq-2} that
\begin{align*}
\begin{aligned} 
\frac{d}{dt}H(\rho^{\eps}\mid M_{\Omega^{\varepsilon}_{\rho^{\eps}}})\leq-\frac{2(d-2)}{e^{2}}\,H(\rho^{\eps}\mid M_{\Omega^{\varepsilon}_{\rho^{\eps}}})
+\frac{d}{dt}\frac{\varepsilon}{\sqrt{|J_{\rho^{\eps}}|^{2}+\varepsilon}}.
\end{aligned}
\end{align*}
Integrating this inequality, we get
\begin{align*}
\begin{aligned} 
H(\rho^{\varepsilon}(t)\mid M_{\Omega^{\varepsilon}_{\rho^{\eps}(t)}})&\leq e^{-C_1t}H(\rho_{0}\mid M_{\Omega^{\eps}_{\rho_0}})
+ e^{-C_1t}\int_0^t e^{C_1s}\,\frac{d}{ds}\frac{\varepsilon}{\sqrt{|J_{\rho^{\eps}(s)}|^{2}+\varepsilon}}\,ds\\
&=e^{-C_1t}H(\rho_{0}\mid M_{\Omega^{\eps}_{\rho_0}})
+ \frac{\varepsilon}{\sqrt{\mid J_{\rho^{\eps}(t)}\mid^{2}+\varepsilon}} - e^{-C_1t} \frac{\varepsilon}{\sqrt{\mid J_{\rho_0}\mid^{2}+\varepsilon}}\\
&\quad -C_1 e^{-C_1t}\int_0^t e^{C_1s}\frac{\varepsilon}{\sqrt{|J_{\rho^{\eps}(s)}|^{2}+\varepsilon}}\,ds\\
&\leq e^{-C_1t}H(\rho_{0}\mid M_{\Omega^{\eps}_{\rho_0}})+\sqrt{\eps}.
\end{aligned}
\end{align*}
%Thus we have
%\[
%H(\rho^{\varepsilon}|M_{\Omega^{\varepsilon}_{\rho^{\eps}}})\leq e^{-C_1t}H(\rho_{0}|M_{\Omega^{\eps}_{\rho_0}})+\frac{\varepsilon (d-1)|J_{\rho^{\eps}}|^2}{C_1(|J_{\rho^{\eps}}|^{2}+\varepsilon)^{3/2}} (1-e^{-C_1t}).
%\]
%Using Lemma \ref{lem-moment} and $|J_{\rho}|\le \int \rho \,d\omega=1$, we have
%\[
%H(\rho^{\varepsilon}|M_{\Omega^{\varepsilon}_{\rho^{\eps}}})\leq e^{-C_1t}H(\rho_{0}|M_{\Omega^{\eps}_{\rho_0}})+ \frac{\varepsilon (d-1)}{C_1(\mid J_{\rho_{0}}\mid^{2}e^{-2(d-2)t}+\varepsilon)^{3/2}} (1-e^{-C_1t}).
%\]
Hence, thanks to the Csiszar-Kullback-Pinsker inequality (see for example \cite[Theorem 1.4]{entropy-inequality}) and the bound
\[
H(\rho_{0}\mid M_{\Omega^{\eps}_{\rho_0}})\le \int_{\bbs^{d-1}}\rho_0\log\rho_0\,d\omega +1-\log C_M,
\]
we have the desired inequality.
\end{proof}

The above estimates immediately imply that our weak solutions of \eqref{main} looks more
and more as a Fisher-von Mises distribution as $t \to \infty$.
\begin{proposition} \label{prop-expo}
Let $\rho$ be a solution of \eqref{main}. Then for all $t>0$,
\[
\|\rho(t)-M_{\Omega_{\rho(t)}}\|_{L^{1}(\mathbb{S}^{d-1})}
 \leq e^{-\frac{2(d-2)}{e^{2}}t}\Big(\int_{\bbs^{d-1}}\rho_0\log\rho_0\,d\omega +1-\log C_M\Big).
\]
\end{proposition}
\begin{proof}
The desired inequality follows by taking $\eps\to 0$ in Lemma \ref{lem-equili} for each $t>0$.
\end{proof}
The above proposition only tells us that our solution $\rho(t)$ resembles to $M_{\Omega_{\rho(t)}}$ for $t \gg 1$,
but it does not say whether the vector $\Omega_{\rho(t)}$ stabilizes to a fixed vector as $t\to \infty$.

To prove this fact, we first use the above result to obtain the uniform positivity of $|J_{\rho}|$ in time, which improves the estimate \eqref{J-0}. The following result ensures that if there is a limit $J_{\infty}$ of  $J_{{\rho}(t)}$ as $t\to\infty$, then $J_{\infty}$ has to be a nonzero vector.
\begin{lemma} \label{lem-positive}
Let $\rho$ be a solution of \eqref{main} with initial data $\rho_{0}\in\mathcal{P}(\bbs^{d-1})$ satisfying \eqref{ini-assume}. Then there exists a positive constant $C(\rho_0)$ only depending on $\rho_{0}$ such that for all $t>0$,
\[
|J_{{\rho}(t)}|>C(\rho_0).
\]
\end{lemma}
\begin{proof}
Using Proposition \ref{prop-expo}, we have
\begin{align*}
\begin{aligned}
\bigg| J_{\rho(t)}-\int_{\bbs^{d-1}}\omega\cdot M_{\Omega_{\rho(t)}} \,d\omega \bigg|
&
= \bigg|\int_{\bbs^{d-1}}\omega\,(\rho(t)-M_{\Omega_{\rho(t)}})\,d\omega\bigg|
\\
&\leq \|\rho(t)-M_{\Omega_{\rho(t)}}\|_{L^{1}(\mathbb{S}^{d-1})}\\
&\leq e^{-\frac{2(d-2)}{e^{2}}t}\Big(\int_{\bbs^{d-1}}\rho_0\log\rho_0\,d\omega +1-\log C_M\Big),
\end{aligned}
\end{align*}  
which yields
\begin{align*}
\begin{aligned}
|J_{\rho(t)}|
&
\geq \bigg|\int_{\bbs^{d-1}}\omega\cdot M_{\Omega_{\rho(t)}} \,d\omega\bigg|-\bigg| J_{\rho(t)}-\int_{\bbs^{d-1}}\omega\cdot M_{\Omega_{\rho(t)}} \,d\omega \bigg|
\\
&\geq \bigg|\int_{\bbs^{d-1}}\omega\cdot M_{\Omega_{\rho(t)}} \,d\omega\bigg|-e^{-\frac{2(d-2)}{e^{2}}t}\Big(\int_{\bbs^{d-1}}\rho_0\log\rho_0\,d\omega +1-\log C_M\Big)\\
&=:R(t).
\end{aligned}
\end{align*}
By \eqref{J-equ}, $C(d):=\big|\int_{\bbs^{d-1}}\omega\cdot M_{\Omega_{\rho}} \,d\omega\big|$ is a positive constant independent of ${\Omega_{\rho}}$, thus 
\[
R(t)\to C(d)\quad\mbox{as}~t\to \infty.
\]
Recalling \eqref{J-0}, this completes the proof.
\end{proof}

\subsection{Proof of Theorem \ref{thm-converge}}
We begin with
\begin{align*}
\begin{aligned}
\frac{d}{dt} J_{\rho}&=\int \omega\, \nabla_{\omega}\cdot\big(\rho \nabla_{\omega}(\log \rho-\nabla_{\omega}(\omega\cdot\Omega_{\rho}))\big) \,d\omega\\
&=-\int \rho\, \nabla_{\omega}\log \rho \,d\omega +\int \rho\,\nabla_{\omega}(\omega\cdot\Omega_{\rho}) \,d\omega.
\end{aligned}
\end{align*}  
If we regard the above terms as functionals on $\rho$, that is,
\begin{align*}
\begin{aligned}
&\mathcal{I}_1(\rho):=\int \rho\, \nabla_{\omega}\log \rho \,d\omega,\\
&\mathcal{I}_2(\rho):=\int \rho\,\nabla_{\omega}(\omega\cdot\Omega_{\rho}) \,d\omega,
\end{aligned}
\end{align*}  
then we see that
\[
\mathcal{I}_1(M_{\Omega_{\rho}})=\mathcal{I}_2(M_{\Omega_{\rho}}).
\]
Also, noticing that $\mathcal{I}_1(\rho)$ and $\mathcal{I}_1(\rho)$ can be written as (see Section \ref{secf:formulas})
\begin{align*}
\begin{aligned}
&\mathcal{I}_1(\rho)=\int \nabla_{\omega} \rho \,d\omega=(d-1)\int \omega\, \rho \,d\omega,\\
&\mathcal{I}_2(\rho)=\int \rho\, \bbp_{\omega^{\perp}}\Omega_{\rho} \,d\omega,
\end{aligned}
\end{align*}
we have
\begin{align*}
\begin{aligned}
|\mathcal{I}_1(\rho)-\mathcal{I}_1(M_{\Omega_{\rho}})|&\le (d-1)\Big| \int \omega \,(\rho-M_{\Omega_{\rho}}) \,d\omega \Big|\\
&\le C \|\rho - M_{\Omega_{\rho}}\|_{L^1(\bbs^{d-1})}, 
\end{aligned}
\end{align*}
and 
\begin{align*}
\begin{aligned}
|\mathcal{I}_2(\rho)-\mathcal{I}_2(M_{\Omega_{\rho}})|&\le C \|\rho - M_{\Omega_{\rho}}\|_{L^1(\bbs^{d-1})}
\end{aligned}
\end{align*}
for some dimensional constant $C$.
Thus, thanks to Proposition \ref{prop-expo}, we have
\begin{align*}
\begin{aligned}
\Big|\frac{d}{dt} J_{\rho}\Big|&=\Big| -\mathcal{I}_1(\rho) + \mathcal{I}_1(M_{\Omega_{\rho}}) -\mathcal{I}_2(M_{\Omega_{\rho}}) + \mathcal{I}_2(\rho) \Big|\\
&\le C \|\rho - M_{\Omega_{\rho}}\|_{L^1(\bbs^{d-1})}.\\
&\le C \,e^{-\frac{2(d-2)}{e^{2}}t}.
%&\to 0\quad \mbox{as}~t\to\infty.
\end{aligned}
\end{align*} 
Together with Lemma \ref{lem-positive}, this implies that there exist a nonzero constant vector $J_{\infty}\in\bbr^d$ such that
\[
 |J_{\rho(t)} -  J_{\infty}| \le C e^{-\frac{2(d-2)}{e^{2}}t}.
\]
Therefore, setting $\Omega_{\infty}=\frac{J_{\infty}}{|J_{\infty}|}$, we have
\[
\|M_{\Omega_{\rho(t)}}-M_{\Omega_{\infty}}\|_{L^1(\bbs^{d-1})}\le C\,|\Omega_{\rho(t)} -\Omega_{\infty}|\le C\, e^{-\frac{2(d-2)}{e^{2}}t},
\]
that combined with Proposition \ref{prop-expo} completes the proof.
\qed

\section{Uniqueness and Stability}
\setcounter{equation}{0}
In this section we present a stability estimate in Wasserstein distance, which provides as a corollary the uniqueness result in Theorem \ref{thm-exist}. First of all, notice that the stability estimate \eqref{eps-uni}  does not imply the stability for \eqref{main}, because of the dependence 
 on $\eps$ in
\eqref{eps-uni}. 

We obtain here a stability estimate for short time when the two initial data $\rho_{0},\bar{\rho}_{0}$ are close to each other as 
\beq\label{close}
W_{2}(\rho_{0},\bar{\rho}_{0}) \le \frac{|J_{\rho_0}|^2}{16}.
\eeq
To get the stability estimate, we use the following lemma on the continuity of the momentum $J_\rho$
with respect to the density $\rho$.
\begin{lemma}\label{lem-mc} 
Let $\rho,\bar\rho\in\mathcal{P}(\bbs^{d-1})$ be any measures satisfying $|J_{\bar{\rho}}|>0$. Then,
\[
 \big| |J_{\bar{\rho}}|-|J_{\rho}| \big|\leq \frac{2\,W_{2}(\bar\rho,{\rho})}{|J_{\bar{\rho}}|} .
\]
\end{lemma}
\begin{proof}
We follow the same arguments in the proof of Proposition \ref{prop-unique}. Let $\varphi_0$ be a $d^2/2$-convex function such that $\exp_{\omega}(\nabla_{\omega}\varphi_0)$ is the optimal map sending $\rho\,d\omega$ onto $\bar\rho\,d\omega$, and
 consider the unique geodesic $r\mapsto\alpha_{r} \,d\omega$ connecting $\alpha_{0}=\rho$ to $\alpha_{1}=\bar\rho$.\\
Similarly for each $r\in[0,1]$, let $\varphi_r$ be a $d^2/2$-convex function such that $\exp_{\omega}(\nabla_{\omega}\varphi_r)$ is the optimal map 
 sending $\alpha_r\,d\omega$ onto $\bar\rho\,d\omega$,
and consider the geodesic $s\mapsto\alpha_{r,s} \,d\omega$ connecting $\alpha_{r,0}=\alpha_{r}$ to $\alpha_{r,1}=\bar\rho$. Notice that $s\mapsto\alpha_{r,s} \,d\omega$ satisfies the continuity equation in the sense of distributions:
\[
\frac{\partial}{\partial s}\bigg|_{s=0}\alpha_{r,s} = -\nabla_{\omega}\cdot(\alpha_{r,s}\nabla_{\omega}\varphi_r ).
\]
Using the same computations as in \eqref{second-1}, we have
\begin{align*}
\begin{aligned}
\frac{\partial}{\partial h}\bigg|_{h=r}\alpha_{h}=\frac{1}{1-r}\frac{\partial}{\partial s}\bigg|_{s=0}\alpha_{r,s} = -\frac{1}{1-r}\nabla_{\omega}\cdot(\alpha_{r}\nabla_{\omega}\varphi_r ),
\end{aligned}
\end{align*}
thus
\begin{align*}
 \begin{aligned}
 \frac{d}{dt}\bigg|_{h=r}| J({\alpha_{h}})|^2 &= 2J({\alpha_{h}})\cdot  \frac{d}{dt}\bigg|_{h=r}J({\alpha_{h}})\\
 &=\frac{2}{1-r}\int_{\bbs^{d-1}}\nabla_{\omega}(\omega\cdot J(\alpha_{r}))\nabla_{\omega}\varphi_{r}\alpha_{r} \,d\omega\\ 
 &\leq\frac{2}{1-r}\sqrt{\int_{S^{d-1}}\mid\nabla_{\omega}(\omega\cdot  J(\alpha_{r}))
 \mid^{2}\alpha_{r} \,d\omega}\sqrt{\int_{S^{d-1}}\mid\nabla_{\omega}
 \varphi_{r}\mid^{2}\alpha_{r} \,d\omega}\\
 &\leq \frac{2}{1-r} W_{2}(\alpha_{r},\bar\rho)\\
 &\leq 2\,W_{2}(\rho,\bar\rho).\\
 \end{aligned}
\end{align*}
Integrating the above inequality from $r=0$ to $r=1$ we get
\[
| J_{\bar{\rho}} |^2- |J_{\rho}|^2 \leq 2\,W_{2}(\rho,\bar\rho).
\]
Similarly applying the above arguments to another $d^2/2$-convex function $\overline{\varphi}_{0}$ sending
$\overline{\rho} \,d\omega$ onto ${\rho} \,d\omega$
we get
\[
|J_{\rho} |^2- |J_{\bar{\rho}}|^2 \leq 2\,W_{2}(\rho,\bar\rho),
\]
hence
\[
\Big| |J_{\bar{\rho}}|-|J_{\rho}| \Big|\leq \frac{2\,W_{2}(\bar\rho,{\rho})}{|J_{\bar{\rho}}|+|J_{\rho}|} \leq \frac{2\,W_{2}(\bar\rho,{\rho})}{|J_{\bar{\rho}}|}.
\]
\end{proof}

\subsection{Proof of Theorem \ref{thm-stable}}
Since $\rho$ solves the continuity equation
\[
\partial_{t}\rho+\nabla_{\omega}\cdot(\rho\,\nabla_{\omega}(\omega\cdot\Omega_{\rho}-\log\rho))=0,
\]
it follows from the Benamou and Brenier formula \cite{Benamou-Brenier} that, for any $t>0$,
\[
W_{2}^2(\rho(t),\rho_0)
\leq t\int_{0}^{t}\int_{\bbs^{d-1}}|\nabla_{\omega}(\omega\cdot\Omega_{\rho(\tau)}-\log\rho(\tau))|^2 \rho(\tau) \,d\omega \,d\tau.
\]
In addition, since 
\begin{align*}
\begin{aligned}
\frac{d}{dt}\mathcal{E}^{0}(\rho)& =-\int_{\bbs^{d-1}}|\nabla_{\omega}(\log\rho-\omega\cdot\Omega_{\rho})|^2 \rho\,d\omega,
\end{aligned}
\end{align*}
we have
\[
W_{2}^2(\rho(t),\rho_0)
\leq t\int_{0}^{t} \Big(-\frac{d}{d\tau}\mathcal{E}^{0}(\rho) \Big) \,d\tau.
\]
Recalling \eqref{relation-1}, we see that
\[
\frac{d}{dt}\mathcal{E}^0(\rho)=\frac{d}{dt}H(\rho\mid M_{\Omega_{\rho}}).
\]
Thus, for all $t>0$,
\[
W_{2}^2(\rho(t),\rho_0)\le t\,\Bigl(H(\rho_0\mid M_{\Omega(\rho_0)})-H(\rho(t)\mid M_{\Omega_{\rho(t)}})\Bigr) 
\le t \,H(\rho_0\mid M_{\Omega(\rho_0)}).
\]
Analogously
\[
W_{2}^2(\bar\rho(t),\bar\rho_0)\le t\,H(\bar\rho_0\mid M_{\Omega(\bar\rho_0)}).
\]
Therefore, setting 
\[
 \delta:=\frac{|J_{\rho_0}|^4}{2^8 \max \{H(\rho_0\mid M_{\Omega(\rho_0)}),H(\bar\rho_0\mid M_{\Omega(\bar\rho_0)})\}}
\]
we have that, for all $t\le \delta$,
\beq\label{instant}
W_{2}(\rho(t),\rho_0)\leq  \frac{|J_{\rho_0}|^2}{16},\qquad W_{2}(\bar\rho(t),\bar\rho_0) \leq \frac{|J_{\rho_0}|^2}{16}.
\eeq
For each time $t\le\delta$, 
%let $\varphi_0$ be a $d^2/2$-convex function such that $\exp_{\omega}(\nabla_{\omega}\varphi_0)$ is the optimal map as
% \[
% \exp_{\omega}(\nabla_{\omega}\varphi_0)_{\#}\rho(t) \,d\omega=\bar\rho(t) \,d\omega.
% \]
we consider the unique geodesic $r\mapsto\alpha_{r} \,d\omega$ connecting $\alpha_{0}=\rho(t)$ to $\alpha_{1}=\bar\rho(t)$.\\
Then, using \eqref{close} and \eqref{instant}, we have
\begin{align*}
\begin{aligned}
W_{2}(\rho_{0},\alpha_{r})&\leq W_{2}(\rho_{0},\rho({t}))+W_{2}(\rho({t}),\alpha_{r})\\
&\leq W_{2}(\rho_{0},\rho({t}))+W_{2}(\rho({t}),\bar{\rho}({t}))\\
&\leq2\,W_{2}(\rho_{0},\rho({t}))+W_{2}(\rho_{0},\bar{\rho}({t}))\\
&\leq2\,W_{2}(\rho_{0},\rho({t}))+W_{2}(\rho_{0},\bar{\rho}_{0})+W_{2}(\bar{\rho}_{0},\bar{\rho}({t}))\\
&\leq\frac{|J_{\rho_0}|^2}{4},
\end{aligned} 
\end{align*}
and applying Lemma \ref{lem-mc} to $\rho_0$ and $\alpha_r$ we get
\[
 \big| |J_{\rho_0}|-|J_{\alpha_r}| \big|\leq \frac{2\,W_{2}(\rho_0,\alpha_r)}{|J_{\rho_0}|}\le \frac{|J_{\rho_0}|}{2},
\]
thus
\beq\label{alpha-r}
|J_{\alpha_r}|\geq  \frac{|J_{\rho_0}|}{2}.
\eeq
We now compute the second derivative of $\mathcal{E}^{0}$ using \eqref{second-1} and \eqref{Hessian}
with $\eps=0$, and thanks to  
\eqref{alpha-r} we have
\[
\frac{d^{2}}{dh^{2}}\bigg|_{h=r}\mathcal{E}^0(\alpha_{h})\ge -\lambda\, W_{2}^{2}(\rho^{\varepsilon}(t),\bar\rho^{\varepsilon}(t)),
\]
where $\lambda:=(1+2/|J_{\rho_0}|)-(d-2)$.\\
Hence,
using the same arguments in the proof of Proposition \ref{prop-unique}, we deduce that
\beq\label{w-stability}
W_{2}(\rho(t),\bar{\rho}(t))\le e^{\lambda t}W_{2}(\rho_0,\bar{\rho}_0)\qquad \forall\,t \in [0,\delta],
\eeq
as desired.
\qed

\subsection{Proof of the uniqueness in Theorem \ref{thm-exist}}
The short time stability estimate \eqref{w-stability} implies the uniqueness of weak solutions to \eqref{main}. Indeed, if $W_{2}(\rho_0,\bar{\rho}_0)=0$, then $W_{2}(\rho(t),\bar{\rho}(t))=0$ for all $t\le\delta$. Thanks to Lemma \ref{lem-positive} and
\[
\frac{d}{dt}H(\rho\mid M_{\Omega_{\rho}})\leq 0,
\]
a continuation argument implies
$W_{2}(\rho(t),\bar{\rho}(t))=0$ for all $t\ge 0$.
\qed

\begin{appendix}
\setcounter{equation}{0}
\section{}
We here present how to compute explicitely the momentum $J_{M_{\Omega}}$ of the Fisher-von Mises distribution
in the case $d=3$.

Let us fix a reference Cartesian coordinate system with $e_3=\Omega$, and then consider the spherical coordinate system $(\theta,\phi)$ associated with the orthonormal basis $(e_1,e_2,\Omega)$. Then a straightforward compution yields
\[
C_M^{-1}=\int_{\bbs^2}e^{\omega\cdot\Omega} \,d\omega=\int_0^{2\pi} d\phi \int_0^{\pi} e^{\cos\theta} \sin\theta \,d\theta =2\pi (e-e^{-1}).
\]
Moreover, since 
\[
\omega=\sin\theta\cos\phi \,e_1 + \sin\theta\sin\phi \,e_2 +\cos\theta \,\Omega,
\]
we have
\[
J_{M_{\Omega}}= \int_{\bbs^2}\omega \,M_{\Omega}(\omega) \,d\omega=C_M\Omega \int_0^{2\pi} d\phi \int_0^{\pi} \cos\theta \,e^{\cos\theta} \sin\theta
\,d\theta =\frac{2e^{-1}}{e-e^{-1}}\, \Omega.
\]
Similarly using the generalized spherical coordinate system on $\bbs^{d-1}$, we have
\beq\label{J-equ}
J_{M_{\Omega}}=\frac{\int_0^{\pi} \cos\theta\, e^{\cos\theta} \sin^{d-2}\theta \,d\theta}{\int_0^{\pi} e^{\cos\theta} \sin^{d-2}\theta \,d\theta} \,\Omega.
\eeq
Notice that $C_{M}$ and $|J_{M_{\Omega}}|$ are constants only depending on dimension $d$, but independent of $\Omega$. 
\end{appendix}

\bibliographystyle{amsplain}
\bibliography{Homogeneous_Vicsek_final}

\end{document}